\documentclass[a4paper,12pt,reqno]{amsart}
\usepackage{amsmath}%
\usepackage{amsfonts}%
\usepackage{amssymb}%
\usepackage{graphicx}
\usepackage{latexsym,euscript} 
\usepackage[latin1]{inputenc}
\newtheorem{theo}{Theorem} 
\newtheorem{Theorem}{Theorem}[section] 
\newtheorem{lemma}[Theorem]{Lemma} 
\newtheorem{prop}[Theorem]{Proposition} 
\newtheorem{rem}[Theorem]{Remark}

\theoremstyle{plain}

\numberwithin{equation}{section}

\newcommand{\finishproof}{\hfill $\Box$ \vspace{5mm}}
\newcommand{\A}{{\mathbb A}}
\newcommand{\C}{{\mathbb C}} 
\newcommand{\Z}{{\mathbb Z}} 
\newcommand{\N}{{\mathbb N}} 
\newcommand{\R}{{\mathbb R}} 
\newcommand{\T}{{\mathbb T}} 
\newcommand{\D}{{\mathbb D}}

\begin{document}
\begin{center}
\centerline{The paper is dedicated to the memory of Professor Nikola\u{\i} Nekhoroshev.}
\end{center}
\title[Effective Stability]{Gevrey Normal Form and Effective Stability of Lagrangian Tori}
\author{Todor Mitev}
\author{Georgi Popov}\thanks{G.P. partially supported by Agence Nationale de la Recherche, Projet 
"RESONANCES": GIP ANR-06-BLAN-0063-03}

\begin{abstract}
A Gevrey symplectic normal form of an analytic and more generally Gevrey smooth  Hamiltonian near a Lagrangian invariant torus with a Diophantine vector of rotation is obtained.  The normal form implies effective stability of the quasi-periodic motion near the torus.

\end{abstract}

\noindent
\maketitle

\section{Introduction}\label{sec:intro}
The aim of this paper is to obtain a Birkhoff Normal Form (shortly BNF) in Gevrey classes of a Gevrey smooth Hamiltonian near a Kronecker torus $\Lambda$ with a Diophantine vector of rotation. Such a  normal form  implies ``effective stability'' of the quasi-periodic motion  near the invariant torus, that is stability in a finite but exponentially long time interval. As in \cite{P2,P4, PT3} it can be used to obtain a microlocal Quantum Birkhoff Normal Form for the Schrödinger operator $P_h = -h^2\Delta + V(x)$ near  $\Lambda$ and to describe the semi-classical behavior of the corresponding eigenvalues (resonances).

 A Kronecker torus of a smooth Hamiltonian $H$ in a symplectic manifold of dimension $2n$ is  a smooth embedded Lagrangian submanifold $\Lambda$,    diffeomorphic to the flat torus $\T^n:=\R^n/2\pi\Z^n$, which is  invariant with respect to the flow $\Phi^t$ of $H$, and such that the restriction of $\Phi^t$  to $\Lambda$ is smoothly conjugated to the linear flow $g^t_\omega(\varphi) := \varphi +t\omega\,  (\mbox{mod}\, 2\pi)$ on $\T^n$ for some $\omega\in\R^n$. Hereafter, we  suppose that $\omega$ satisfies the usual  Diophantine condition \eqref{eq:sdc}. Then there is a symplectic mapping $\chi$  from a neighborhood of the zero section $\T^n_0:=\T^n\times \{0\}$ to  a neighborhood of $\Lambda$ in $X$ sending $\T^n_0$ to $\Lambda$ and such that the Hamiltonian  $H_0:=\chi^\ast H$ becomes  $H_0(\varphi,I) = H^0(I) + R^0(\varphi,I)$, where $\nabla H^0(0)=\omega$, and the Taylor series of $R^0$ at $I=0$ vanishes (cf. \cite{La}, Proposition 9.13). In particular,  
$\T^n_0$ is an invariant torus of $H_0$, the restriction of the flow of $H_0$ to $\T^n_0$ is given by $g^t_\omega(\varphi) = \varphi +t\omega\,  (\mbox{mod}\, 2\pi)$,  and for any $\alpha,\, \beta\in \N$, and any $N\ge 1$, we have 
$\partial_\varphi^\alpha \partial_I^\beta R^0(\varphi,I)= O_{\alpha,\beta,N}(|I|^N)$. Our aim is to replace these polynomial estimates  by exponential estimates of $\partial_\varphi^\alpha \partial_I^\beta R^0$ of the form $O_{\alpha,\beta}(\mbox{exp}(-c/|I|^a))$, $a>0$, $c>0$,  in the case when both the Hamiltonian $H$  and the Kronecker torus are Gevrey smooth. In the case when the Hamiltonian and the torus are analytic  a similar BNF has been obtained by Morbidelli and  Giorgilli. They have proved as well  effective stability of the action  near  analytic KAM tori and even a super exponential stability of the action  \cite{G-M1, M-G1, M-G2} for convex Hamiltonians using Nekhoroshev's theory. 
A simultaneous normal form for a family of Gevrey KAM tori has been obtained in \cite{P1,P3}. 

The existence of large family of Kronecker  tori of Diophantine vectors of rotation is given by the classical KAM theorem in the case of real-analytic Hamiltonians satisfying the Kolmogorov non-degeneracy condition and for Gevrey smooth Hamiltonians it has been proved by one of the authors in \cite{P3}. 
Similar results for analytic (Gevrey-smooth) Hamiltonians satisfying  the Rüssmann non-degeneracy conditions have been obtained in \cite{X-Y1}. 

\section{Main results}\label{sec:main}
Let $X$ be a bounded domain in ${\R}^n$. Fix $\rho \ge 1$ and a positive constant $L$, and denote by 
${\mathcal G}^{\rho}_{L}(X)$ 
the set of all $C^\infty$-smooth  
real-valued functions   $H$ in $X$  such that 
\begin{equation}
\|H\|_{L} := 
\sup_{\alpha\in{\N}^n}\, 
\sup_{x\in X }\, 
\left(|\partial_x^\alpha 
H(x)|\,  L^{-|\alpha|} 
\alpha !^{-\rho} \right)  < \infty \, ,	
                                    \label{eq:gevrey-1}
\end{equation}
where $\N$ is the set of non negative integers, $\alpha ! =
\alpha_1! \cdots \alpha_n!$ and  $|\alpha| = \alpha_1 + \cdots + \alpha_n$ is the ``length'' of  $\alpha=(\alpha_1,\ldots,\alpha_n)
\in {\N}^n$. A function $H$ is said to be 
${\mathcal G}^{\rho}$-smooth  on $X$ if it satisfies (\ref{eq:gevrey-1}) with some $L>0$. 
In the same way, using local coordinates,  we define ${\mathcal G}^{\rho}$-smooth functions on a  ${\mathcal G}^{\rho}$-smooth manifold $X$ of dimension $n$. Note that the ${\mathcal G}^{1}$-smooth functions in a bounded domain (real-analytic manifold) $X$   are just the analytic functions in $X$. On the other hand,  the class of ${\mathcal G}^{\rho}$-smooth function is not quasi-analytic for $\rho>1$;  there exist  functions of a compact support which are ${\mathcal G}^{\rho}$-smooth. 
For more properties of  Gevrey smooth functions  we refer to \cite{Kom, L-M} and \cite[Appendix]{P3}, where the implicit function theorem and the composition of Gevrey functions is discussed. 

When dealing with the KAM theory in Gevrey classes one looses Gevrey regularity in  frequencies, and there naturally arise anisotropic Gevrey classes. They  are defined
as follows. 
Let $\rho,\mu \ge 1$ and $L_1, L_2$
be positive constants.  
Given a bounded domain $D \subset {\R}^n$, we consider  $\A:=\T^n\times D$ provided with the canonical symplectic structure, and 
denote by 
${\mathcal G}^{\rho,\mu}_{L_1,L_2}(\A)$ 
the set of all $C^\infty$-smooth 
real valued Hamiltonians  $H$ in $\A$  such that 
\begin{equation}
\|H\|_{L_1,L_2} := 
\sup_{\alpha,\beta\in \N^n}
\sup_{(\theta, I)\in \A } 
\left(|\partial_\theta^\alpha \partial_I^\beta  
H(\theta, I)|\,  L_1^{-|\alpha|} L_2^{-|\beta|}
\alpha !^{-\rho} \beta !^{-\mu}\right) < \infty .	
                                    \label{eq:gevrey-2}
\end{equation}
 A Hamiltonian $H$ in $\A$ is said to be 
${\mathcal G}^{\rho,\mu}$-smooth if it belongs to ${\mathcal G}^{\rho,\mu}_{L_1,L_2}(\A)$ for some positive constants $L_1, L_2$. The numbers  $\rho\ge 1$ and $\mu\ge 1$ in \eqref{eq:gevrey-2} are called Gevrey exponents and the positive  constants $L_1$ and $L_2$ are called Gevrey constants. We say that the two pairs of Gevrey constants $L_1$, $L_2$ and   $\tilde L_1$, $\tilde L_2$ are {\em equivalent} if there is $c_1(n,\rho,\mu)>0$ and $c_2(n,\rho,\mu)>0$ such that $\tilde L_1 = c_1(n,\rho,\mu) L_1$ and $\tilde L_2=c_2(n,\rho,\mu)L_2 $. 

Let $\rho\ge 1$ and  let $X$ be a  ${\mathcal G}^{\rho}$-smooth symplectic manifold of dimension $2n$.  
 Let $H$ be a ${\mathcal G}^{\rho}$-smooth Hamiltonian in $X$. 
A ${\mathcal G}^{\rho}$-smooth Kronecker torus  of $H$ of a vector of rotation $\omega\in\R^n$ is given by a 
${\mathcal G}^{\rho}$-smooth embedding $f:\T^n\to X$, such that $\Lambda = f(\T^n)$ is a Lagrangian submanifold of $X$ which is   invariant with respect to the Hamiltonian flow $\Phi^t$ of $H$ and  $\Phi^t\circ f = f\circ g^t_\omega$ for all $t\in\R$, i.e. the  
 diagram
\begin{equation}
\label{eq:diagram} 
\begin{array}{cccl} 
\T^n&\stackrel{g^t_\omega}{\longrightarrow}&\T^n\cr
\downarrow\lefteqn{f}& &\downarrow\lefteqn{f} \cr
\Lambda&\stackrel{\Phi^t}{\longrightarrow}&\Lambda&  
\end{array} 
\end{equation}
is commutative for any $t\in\R$. Recall that $g^t_\omega(\varphi) =\varphi + t\omega\ (\mbox{mod}\,  2\pi)$. 
We will suppose that $\omega$ is $(\kappa,\tau)$-Diophantine for some $\kappa >0$
and $\tau > n-1$, which means the following:
\begin{equation}
\mbox{For any}\  0\neq k\in \Z^{n}\, , \quad                    
|\langle\omega,k\rangle |\ \ge \ \kappa\, 
|k|^{- \tau}\, ,
                       \label{eq:sdc}                                        
\end{equation}
where $|k|=\sum_{j=1}^n |k_j|$. 
Note that if  $X$ is exact symplectic and $\Lambda\subset   X$ is an  embedded submanifold 
  satisfying (\ref{eq:diagram}) with a Diophantine vector $\omega$ then $\Lambda$ is
Lagrangian (see \cite{Her}, Sect. I.3.2).  The existence of such tori in  $\A:=\T^n\times D$ 
with vectors of rotation $\omega$ satisfying \eqref{eq:sdc} is provided by the KAM theorem.
It follows from \cite[Theorems 1.1 and 3.12]{P3} and \cite{P1},  that if $H\in {\mathcal G}^{\rho}(\A)$ is a ``small'' (in terms of $\kappa$) real-valued  perturbation of a completely integrable Hamiltonian satisfying the Kolmogorov non-degeneracy conditions, then there is a Cantor set $\Omega_\kappa\in\R^n$ of frequencies satisfying \eqref{eq:sdc} and of a positive Lebesgue measure such that  for any $\omega\in \Omega_\kappa$ there is a ${\mathcal G}^{\rho}$-smooth  Kronecker torus $\Lambda_\omega$ with frequency $\omega$. In the analytic case ($\rho=1$) this follows from the classical KAM theorem. Moreover, the family $\Lambda_\omega$,  $\omega\in \Omega_\kappa$, is ${\mathcal G}^{\mu}$-smooth in Whitney sense, where $\mu=\rho(\tau+1)+1$ when $\rho>1$ and $\mu$ could be  any number greater than $\tau+2$ when $\rho=1$ (see 
\cite{P1, P3, W}). 
The main result in this paper is concerned with a Gevrey smooth  Birkhoff Normal Form of $H$ near any Kronecker torus with a Diophantine frequency. 
\begin{theo}\label{theo:main}
Let $\omega\in\R^n$ satisfy the $(\kappa,\tau)$-Diophantine condition (\ref{eq:sdc}) with some $\kappa>0$ and $\tau>n-1$.  
Fix $\rho \ge  1$ and  set $\mu = \rho(\tau+1)+1$. 
Let $H\in {\mathcal G}^{\rho}(X,\R)$  be a real-valued Hamiltonian and 
let $\Lambda$ be a ${\mathcal G}^{\rho}$-smooth Kronecker torus of $H$ of a vector of rotation $\omega$.  
Then there is a neighborhood  $D$  of $0$ in $\R^n$ and a  symplectic mapping 
$\chi\in {\mathcal G}^{\rho,\mu}(\A, X)$, where $\A=\T^n\times D$, 
such that $\chi(\T^n_0)= \Lambda$, and 
\[
\left\{
\begin{array}{lcrr}
H(\chi(\varphi,I)) = H^0(I)+ R^0(\varphi,I) ,\ \mbox{where}\ 
 H^0\in {\mathcal G}^{\mu}(D ),\   R^0\in {\mathcal G}^{\rho,\mu}(\A ), 
 \\[0.3cm]   
\mbox{and}\ \partial_I^\alpha R^0(\varphi, 0) = 0\ \mbox{for any $\varphi\in\T^n$ and 
 $\alpha\in \N^n$}\, .
 \end{array}
 \right.
 \] 
\end{theo}

In the analytic case ($\rho=1$), a similar BNF near an elliptic equilibrium point of the Hamiltonian  has been obtained by Giorgilli, Delshams, Fontich,  Galgani and Sim\'{o}  in \cite{G-D-F-G-S}. Moreover,  effective stability of the action, that is stability of the action in a finite but exponentially long time interval has been proved in \cite{G-D-F-G-S}. 
Effective stability near an analytic KAM torus has been investigated by  Morbidelli and  Giorgilli in \cite{M-G1},  \cite{M-G2} and  \cite{G-M1}. Combining it with the Nekhoroshev theorem they obtained  a super-exponential effective stability of the action near the torus. 
The Nekhoroshev theory for Gevrey smooth Hamiltonians has been developed by J.-P. Marco and D. Sauzin \cite{M-S}.
As it was mentioned above,  if $H\in {\mathcal G}^{\rho}(\A)$ is a ``small'' (in terms of $\kappa$) real-valued  perturbation of a completely integrable ${\mathcal G}^{\rho}$-smooth Hamiltonian satisfying the Kolmogorov non-degeneracy conditions, then there is a Cantor set $\Omega_\kappa\subset\R^n$ of frequencies satisfying \eqref{eq:sdc} of  positive Lebesgue measure such that  for any $\omega\in \Omega_\kappa$ there is a ${\mathcal G}^{\rho}$-smooth  Kronecker torus $\Lambda_\omega$ with frequency $\omega$.  The family $\Lambda_\omega$,  $\omega\in \Omega_\kappa$, is ${\mathcal G}^{\mu}$-smooth in Whitney sense, where $\mu=\rho(\tau+1)+1$ if $\rho>1$ and $\mu>\tau+2$ if $\rho=1$ (see 
\cite{P1, P3, W}). This implies a 
 simultaneous ${\mathcal G}^{\rho, \mu}$-smooth BNF of the corresponding Hamiltonian at  a family of KAM tori 
 $\Lambda_\omega$, $\omega\in \widetilde \Omega_\kappa$, where 
 $\widetilde \Omega_\kappa\subset \Omega_\kappa $ is the set of points of positive Lebesgue density in $ \Omega_\kappa $  \cite[Corollary 1.2]{P3}. Normal forms for reversible analytic vector fields with an exponentially small error term have been obtained by Iooss and   Lombardi \cite{I-L1, I-L2}.
 
 Here we obtain a BNF of any single  ${\mathcal G}^{\rho}$-smooth  Kronecker torus $\Lambda$ of the Hamiltonian.
This normal form implies effective stability not only of the action but of the {\em quasi-periodic motion} near $\Lambda$ as well (cf. \cite[Corollary 1.3]{P3}).
Moreover, our method allows us to keep track on the corresponding Gevrey constants. 
 In the case of KAM tori \cite{P1, P3} this yields an uniform bound on the corresponding Gevrey constants with respect to  $\omega\in \Omega_\kappa$. It could be applied as in \cite{M-G1},  \cite{M-G2} and  \cite{G-M1} to obtain a super-exponential effective stability of the action near the torus in the case of convex Hamiltonians using the Nekhoroshev theory for Gevrey Hamiltonians developed  by J.-P. Marco and D. Sauzin \cite{M-S}.
It seems that this method could be applied to obtain a Gevrey normal form in the case of elliptic tori and near an elliptic equilibrium point of Gevrey smooth Hamiltonians as well as in the case of hyperbolic tori and reversible systems.  

The method we use relies on an explicit construction of the generating function of the canonical transformation putting the Hamiltonian in a normal form which allows us to obtain an explicit form of the corresponding homological equation (see Sect. \ref{sec:homological}). It is different from those used in \cite{G-D-F-G-S} and \cite{M-G1} which is based on the formalism of the Lie transform. 

It is an interesting question if the exponent $\mu=\rho(\tau +1)+1$ is optimal. As it was mentioned above the same exponent appears in the KAM theorem in Gevrey classes when $\rho>1$ and our exponent $\mu$  is smaller when $\rho=1$, in particular we obtain the same exponent as in the  simultaneous BNF of the family of KAM tori 
 $\Lambda_\omega$, $\omega\in \widetilde \Omega_\kappa$ in  \cite[Corollary 1.2]{P3} when $\rho>1$. In the analytic case ($\rho=1$) there is an heuristic argument of  Morbidelli and  Giorgilli  \cite[§3. Discussion]{M-G2} showing that 
 $\mu=\tau +2$ should be optimal.

 Theorem \ref{theo:main} can be used as in \cite{P2,P4, PT3} to obtain a microlocal Quantum Birkhoff Normal Form in Gevrey classes for the Schrödinger operator $P_h = -h^2\Delta + V(x)$ near a Gevrey smooth Kronecker torus $\Lambda$ of the Hamiltonian $H(x,\xi) = \|\xi\|^2 + V(x)$.  

\section{Birkhoff Normal Form in Gevrey classes and Effective Stability}\label{sec:BNF} 

We are going to reduce the problem to the case of a Gevrey smooth (real-analytic) Hamiltonian in $\A = \T^n\times D$ having a Kronecker torus $\T^n_0=\T^n\times \{0\}$, where $D$ is a connected neighborhood of $0$ in $\R^n$ and $\A$ is provided with the canonical symplectic two-form.  By a result of Weinstein  there  is a symplectic transformation 
$\chi_0:\A \to X$ such that $\chi_0(\T^n_0)=\Lambda$ and $\chi_0\circ \imath = f$, where $\imath(\theta)=(\theta,0)\in\T^n_0$ for any $\theta\in\T^n$. 
To construct $\chi_0$ we first  find  a tubular neighborhood $U$ of $\Lambda$ in $T^\ast \Lambda$ and a ${\mathcal G}^{\rho}$-smooth    symplectic transformation  $F: U\to X$ which maps  the zero section of $\Lambda$ in $T^\ast \Lambda$ to $\Lambda$. If $\rho>0$ one just follows the proof of  Weinstein. 
In the real-analytic case ($\rho=1$), we first take a $C^\infty$-smooth symplectic map $F_0$ with this property, which exists by the Weinstein theorem, next we approximate  it with a real-analytic one, and then we use a deformation argument of Moser to get $F$. Set $\tilde f = F^{-1}\circ f$. 
Arguing as in the proof of Proposition 9.13 \cite{La},   we obtain  a
${\mathcal G}^{\rho}$-smooth symplectic mapping $\chi_1$ from a 
 bounded neighborhood  $\A = \T^n\times D$ of the torus $\T^n_0$ in $T^\ast \T^n$ to a tubular  neighborhood of the zero section of $\Lambda$ in $T^\ast\Lambda$ such that  $\chi_1\circ\imath=\tilde f$, and we set $\chi_0= F\circ \chi_1$.  In particular,  $\chi_0(\T^n_0)=\Lambda$. 
 Moreover, the pull-back of the Hamiltonian vector field to $\A$ is globally Hamiltonian and we denote by $H\in {\mathcal G}^{\rho}(\A, \R)$  its  Hamiltonian  in $\A$. It follows from \eqref{eq:diagram} that the restriction of the flow of the Hamiltonian vector field of $H$ to $\T^n_0$ is just $g^t_\omega$. Moreover, 
 $H(\theta,0)$ is  constant since the flow is transitive in $\T^n_0$, and we take it to be zero. Hence, 
\begin{equation}
\label{eq:H}
H(\theta,r) = \langle\omega, r\rangle + \widetilde H(\theta,r),\   \widetilde H \in {\mathcal G}^{\rho,\rho}(\A), \  \widetilde H(\theta,r)=O(|r|^2) .
\end{equation}
Denote by $\Gamma(t)$, $t>0$, the Gamma function \eqref{eq:GammaFunction}.   Using  Remark \ref{rem:gamma}, we write the corresponding Gevrey estimates as follows 
\begin{equation}
\label{eq:gevrey-coefficients0}
|\partial^\beta_\theta \partial^\alpha_r\widetilde H(\theta,r)| \le  L_0  L_1^{|\beta|}  L_2^{|\alpha|-1}\,   \alpha !\,  \Gamma(\rho|\beta|+1)\Gamma( (\rho-1)|\alpha| +1)\,  
\end{equation}
for any $(\theta,r)\in \A$ and $\alpha,\beta\in\N^n$, where $L_0$, $L_1$ and $L_2$ are positive constants, and we suppose that  $L_0\ge 1$, 
$L_1\ge 1$ and $L_2\ge 1$.  

A smooth function $g(\theta,I)$ in  $\A' = \T^n\times D'$  is said to be a generating function of a canonical transformation $\chi:\A'\to\A$ if 
\begin{equation}
\label{eq:canonical}
\begin{array}{lcrr}
{\rm graph}\, \chi := \{(\chi(\varphi,I);(\varphi,I)):\, (\varphi,I)\in \A'\}\\[0.3cm]
= \displaystyle \left\{ \left(\theta, I+\frac{\partial g}{\partial \theta}(\theta, I); 
\theta +\frac{\partial g}{\partial I}(\theta, I), I\right) \right\}\, . 
\end{array}
\end{equation}
Without loss of generality we can suppose that $\kappa\le 1$ in \eqref{eq:sdc}. 
Theorem \ref{theo:main} follows from the following 

\begin{theo}\label{theo:gevrey}
Let  $\rho \ge   1$ and   $H\in {\mathcal G}^{\rho,\rho}(\A, \R)$. Suppose that $H$ satisfies  \eqref{eq:H} and \eqref{eq:gevrey-coefficients0}, where $\omega\in\R^n$ is $(\kappa,\tau)$-Diophantine and $0< \kappa \le 1$ and $\tau>n-1$. Set  $\mu = \rho(\tau+1)+1$.
 Then there is a neighborhood $D'$ of $0$ in $\R^n$ and a function $g\in {\mathcal G}_{C_1,C_2}^{\rho,\mu}(\A', \R)$, $g(\theta,I) = O(|I|^2)$ in $\A' =\T^n\times D'$, generating a canonical transformation $\chi\in {\mathcal G}^{\rho,\mu}(\A', \A)$,  such that 
\begin{equation}\label{eq:normal-form1}
\left\{
\begin{array}{lcrr}
H(\chi(\varphi, I)) = H^0(I) + R^0(\varphi, I),\\[0.3cm]
\mbox{where}\ H^0\in {\mathcal G}^{\mu}_{C_2}(D',\R),\ R^0\in {\mathcal G}^{\rho,\mu}(\A',\R),\\[0.3cm]
\mbox{and}\  \partial^\alpha_I R^0(\theta, 0) = 0\ \mbox{for any}\ \alpha \in \N^n. 
\end{array}
\right.
\end{equation}
Moreover, the Gevrey constants $C_1$ and  $C_2$ are equivalent to $L_1$ and $\frac{1}{\kappa} L_0L_1^{\tau+n +4} L_2$ respectively, i.e. 
\begin{equation}
\label{eq:gevrey-constants1}
C_1 =c_1(\rho,\tau,n) L_1 \quad \mbox{and} \quad C_2 =c_2(\rho,\tau,n)\frac{1}{\kappa} L_0L_1^{\tau+n +4} L_2\, ,
\end{equation} 
where $c_1$ and $c_2$ are  positive constant depending only on $\rho$, $\tau$ and $n$, while $\kappa$ is the constant in \eqref{eq:sdc}.  
\end{theo}
\noindent
\begin{rem}\label{rm:gevrey-constants}
We have $\chi\in {\mathcal G}_{  C_1,C_2}^{\rho,\mu}(\A', \A)$ and 
$R^0\in {\mathcal G}^{\rho,\mu}_{C_1,C_2}(\A',\R)$, 
 where the Gevrey constants $C_1$ and  $C_2$ are equivalent to $L_1^2$ and $\frac{1}{\kappa} L_0L_1^{\tau+n +6} L_2$ respectively, i.e. 
 \begin{equation}
\label{eq:gevrey-constants2}
C_1 =c_1(\rho,\tau,n) L_1^2 \quad \mbox{and} \quad C_2 =c_2(\rho,\tau,n)\frac{1}{\kappa} L_0L_1^{\tau +n+6} L_2\, ,
\end{equation} 
\end{rem}
Theorem \ref{theo:gevrey} and Remark \ref{rm:gevrey-constants} will be proved in Sect. \ref{sec:Gevrey}. 

By the 
Taylor formula of order $m$ applied to $R^0(\varphi, I)$ at  $I=0$  we obtain
for any $\alpha,\beta \in {\N}^n$, $m\in {\N}$, and $(\varphi,I)\in {\T}^n\times D'$ the estimate
$$
|\partial_\varphi^\alpha\partial_{I}^{\beta} R^0(\varphi, I)|\ \leq\ 
A\,  C_1^{|\alpha|} 
C_2^{|\beta|+m}\, 
\alpha !\,^\rho  \beta !\, ^{\mu}\, m !\, ^{\mu-1}|I|^m\, ,   
$$ 
where $A>0$ and the positive constants $C_1$ and $C_2$ are as in (\ref{eq:gevrey-constants2}). 
Using Stirling's formula 
we minimize the right-hand side  with respect to 
$m\in {\N}$. An optimal choice for $m$ will be 
\[
m\sim \left( C_2|I|\right)^{\, -\frac{1}{\rho(\tau + 1)}},
\]
 which leads to  
\begin{equation}
\begin{array}{rcll}
\displaystyle |\partial_\varphi^\alpha\partial_{I}^{\beta} R^0(\varphi, I)| &\leq& 
 A\,  C_1^{|\alpha|} C_2 ^{|\beta|}\, 
\alpha !\,^\rho  \beta !\, ^{\mu-1}\,  \\[0.5cm]
&\times&  \displaystyle 
\exp\left(- \left( C_2 
|I|\right)^{\, -\frac{1}{\rho(\tau + 1)}}\right)  
                            \label{eq:exponential-estimate}
\end{array}
\end{equation}
 for any $\alpha,\beta \in {\N}^n$ uniformly with respect to 
$(\varphi,I)\in {\T}^n\times D'$, where
$C_1$ and  
$C_2$ are of the form  \eqref{eq:gevrey-constants2}.  
This estimate yields    
effective stability of the quasi-periodic motion near 
the invariant  tori as in \cite[Corollary 1.3]{P3}).

\subsection{Idea of the Proof of Theorem \ref{theo:gevrey}}\label{susec:idea}
Expanding $ \widetilde H(\theta,r)$ in Taylor series with respect to $r$ at $r=0$ we obtain 
\begin{equation}
\label{eq:hamiltonian}
H(\theta,r)\,  \sim\,  \langle\omega, r\rangle + \sum_{m=2}^{\infty}\, H_m(\theta,r)\, ,\quad H_m(\theta,r) = \sum_{|\alpha|=m}b_\alpha(\theta)r^\alpha\, .
\end{equation}
It follows from \eqref{eq:gevrey-coefficients0} that  the coefficients $b_\alpha$ satisfy the following Gevrey type estimates
\begin{equation}
\begin{array}{lcrr}
\label{eq:gevrey-coefficients}
|\partial^\beta b_\alpha(\theta)| = (\alpha !)^{-1}|\partial_\theta^\beta\partial_r^\alpha   \widetilde H(\theta,0)|\\[0.3cm]
\le  L_0  L_1^{|\beta|}  L_2^{|\alpha|-1}\,   \Gamma(\rho|\beta|+1)\Gamma( (\rho-1)|\alpha| +1) ,
\end{array}
\end{equation}
for any $\theta\in \T^n$ and any multi-indices $\alpha,\beta\in \N^n$, $|\alpha|\ge 2$.

We are looking for a function $g\in {\mathcal G}^{\rho,\mu}(\A')$, where  $\A' = \T^n\times D'$ and $D'\subset \R^n$ is a neighborhood of $0$, 
such that $g(\theta,0)=0$,  $\nabla_I g(\theta,0)=0$,   and 
\begin{equation}\label{eq:normal-form}
\left\{
\begin{array}{lcrr}
H(\theta, I + \nabla_\theta g (\theta, I)) = H^0(I) + R(\theta, I)\, ,\\[0.3cm]
\mbox{where}\ H^0\in {\mathcal G}^{\mu}(D',\R)\, ,\ R\in {\mathcal G}^{\rho,\mu}(\A',\R) \, ,\\[0.3cm]
\mbox{and}\  \partial^\alpha_I R(\theta, 0) = 0\ \mbox{for any}\ \alpha \in \N^n\, . 
\end{array}
\right.
\end{equation}
If such a function $g$ exists, and if $D'$ is sufficiently small, 
we get by means of the implicit function theorem in anisotropic Gevrey classes \cite{Kom}, \cite[Proposition A.2]{P2},  
a function $\theta(\varphi,I)$ in ${\mathcal G}^{\rho,\mu}(\A')$ which solves the equation 
\[
\varphi = \theta + \nabla_I g(\theta, I)
\]
with respect to $\theta\in \T^n$, and we denote by $\chi$ the canonical transformation  defined by $g$ by means of \eqref{eq:canonical}. 
Hence,  
\[
(H\circ\chi)(\varphi, I) = H^0(I) + R(\theta(\varphi,I),I).
\]
Setting $R^0(\varphi, I) = R(\theta(\varphi,I),I)$ 	we obtain $R^0\in {\mathcal G}^{\rho,\mu}(\A')$ by the theorem of  composition in anisotropic Gevrey classes \cite[Proposition A.4]{P2}, as well as the identities 
\[
\partial^\alpha_I R^0(\theta, 0) = 0\, . 
\]
for any $\alpha \in \N^n$ and $\theta\in\T^n$.  

Theorem \ref{theo:gevrey} follows from the following 

\begin{prop}\label{prop:gevrey}
Let  $\rho \ge   1$, $\tau>n-1$, and    $\mu = \rho(\tau+1)+1$. Suppose that the Hamiltonian $H\in {\mathcal G}^{\rho,\rho}(\A, \R)$ satisfies  \eqref{eq:H} and \eqref{eq:gevrey-coefficients0}, where $\omega$ satisfies \eqref{eq:sdc}. 
 Then there is a neighborhood $D'$ of $0$ in $\R^n$ and a function $g\in {\mathcal G}_{C_1,C_2}^{\rho,\mu}(\A', \R)$, $g(\theta,I) = O(|I|^2)$ in $\A' =\T^n\times D'$,  such that 
\begin{equation}\label{eq:normal-form2}
\left\{
\begin{array}{lcrr}
H(\theta, I +  \nabla_\theta g (\theta,I)) = H^0(I) + R(\theta, I)\, ,\\[0.3cm]
\mbox{where}\ H^0\in {\mathcal G}^{\mu}_{C_2}(D',\R)\, ,\ R\in {\mathcal G}^{\rho,\mu}_{ C_1, C_2}(\A',\R) \, ,\\[0.3cm]
\mbox{and}\  \partial^\alpha_I R(\theta, 0) = 0\ \mbox{for any}\ \alpha \in \N^n\, ,  
\end{array}
\right.
\end{equation}
where
$C_1$ and  
$C_2$ are given by \eqref{eq:gevrey-constants1}.
\end{prop}

\section{Weighted  Wiener norms}\label{sec:Wiener}\ 
To obtain sharp estimates in Gevrey classes  
we will use weighted Wiener norms. These norms are well adapted to solve the so called homological equation and they provide a sharp estimate for the product of two functions. 
Given $u\in C(\T^n)$, we denote by $u_k$, $k\in\Z^n$, the corresponding Fourier coefficients, and  by  
$$\langle u\rangle:=u_0= (2\pi)^{-n}\int_{\T^{n}}u(\varphi)d\varphi $$
 the mean value of $u$ on $\T^n$. 
For any $s\in\R_+:=[0,+\infty)$ we define the corresponding weighted  Wiener norm of $u$ by
\[
S_s(u)\, :=\ \sum_{k\in\Z^n}\, (1+ |k|)^s|u_k|\, ,
\]
where  $|k|=|k_1|+\cdots+|k_n|$, $k=(k_1,\ldots,k_n)\in \Z^n$. 
 The weighted Wiener space 
${\mathcal A}^s(\T^{n})$, $s\ge 0$, is defined as the Banach space of all $u\in C(\T^n)$ such that $S_s(u) <\infty$ equipped with the norm $S_s$.  The space ${\mathcal A}^s(\T^{n-1})$ is  a Banach algebra, if $u,v \in {\mathcal A}^s(\T^{n})$ then $S_s(uv)\le S_s(u)S_s(v)$. 
Moreover, 
the  following relations between Wiener spaces and Hölder spaces hold
\[
C^q(\T^n) \hookrightarrow {\mathcal A}^s(\T^n)\hookrightarrow C^s(\T^n)\, ,
\]
for any $s\ge 0$ and $q>s+n/2$, and the corresponding inclusion maps are continuous. The first relation  is a special case of a theorem of Bernstein ($n=1$) and its generalizations for $n\ge 2$ \cite[Chap. 3, § 3.2]{kn:AAP}. For more properties of these spaces see \cite{PT3}. 

Weighted  Wiener spaces are perfectly adapted for solving the homological equation  
\begin{equation}
{\mathcal L}_\omega u(\varphi) = f(\varphi)\, 
\label{eq:homological1}
\end{equation} 
where ${\mathcal L}_\omega:= \langle\omega, \frac{\partial }{\partial \theta}\rangle$.  
 We have the following  
\begin{lemma}\label{lemma:homological} Let $\omega$ satisfy the $(\kappa,\tau)$-Diophantine condition (\ref{eq:sdc}) and let $s\ge 0$.  Then for any  $f\in {\mathcal A}^{s+\tau}(\T^n)$ such that $\langle f\rangle =0$
the homological equation  
\[
{\mathcal L}_\omega u=f\, ,\quad    \langle u\rangle =0 ,
\]
has an unique   solution $u\in {\mathcal A}^{s}(\T^n)$,  and it satisfies the estimate 
\[
S_{s}(u)\, \le\,  \frac{1}{\kappa}\, S_{s+\tau}(f)\, . 
\]
\end{lemma}
\noindent
{\em Proof}. 
Comparing the  Fourier coefficients $u_k$
and  $f_k$, $k\in \Z^{n}$, of $u$ and $f$ respectively,  we get 
\[
u_k\ =\ \frac{f_k}{i \langle
k,\omega\rangle }\,  , \ k\neq 0\, ,
\] 
and set $u_0 = 0$. Then using 
(\ref{eq:sdc}) we obtain
\[
|u_k| \le  \frac{1}{\kappa} |k|^\tau |f_k|\le  \frac{1}{\kappa} (1+|k|)^\tau |f_k|\,  , \ k\neq 0\, .
\] 
Since  $f_0=\langle f\rangle = 0$,  taking the 
sum with respect to $k\neq 0$ we get the function $u$ and the corresponding estimate of $S_{s}(u)$. 
In this way we obtain an unique solution $u$
of (\ref{eq:homological1}) 
normalized by $\langle u\rangle =0 $.
 \finishproof \\
In what follows we shall need a sharp estimate of the weighted Wiener norm of the product $uv$ of two functions  $u,v\in {\mathcal A}^s(\T^{n})$.  
Let $[s] \in\Z$ be the integer part of $s\in \R$ and denote by $\{s\}=s - [s]\in [0,1)$ its fractional part.
\begin{lemma}\label{lemma:product}
 For any  $s\in \R_+$  and $u,v\in {\mathcal A}^s(\T^{n})$ we have
\[
S_s(uv)\le 2\sum_{m=0}^{[s]} \begin{pmatrix}[s]\cr m\end{pmatrix}
\big[S_{s-m}(u)S_m(v) + S_{s-m}(v)S_m(u)\big]\, .
\]
\end{lemma}
\noindent
{\em Proof}. For any  $k\in\Z^n$ we set $\langle k\rangle := 1+|k|$. Obviously, 
$\langle k\rangle< \langle l\rangle + \langle k-l\rangle$ for any  $k,l\in\Z^n$, and   we obtain 
\[
\begin{array}{lcrr}
\displaystyle  \langle k\rangle^s\left|(uv)_k\right|\ \le\ 
\sum_{l\in\Z^n}\big(\langle l\rangle+ \langle k-l\rangle\big)^{[s]+\{s\}}|u_l||v_{k-l}| 
 \\[0.3cm]
\displaystyle \le \sum_{l\in\Z^n}\sum_{m=0}^{[s]} \begin{pmatrix}[s]\cr m\end{pmatrix} 
\big(\langle l\rangle+ \langle k-l\rangle\big)^{\{s\}}\langle l\rangle^{m}|u_l|\langle k-l\rangle^{[s]-m}|v_{k-l}|
\\[0.3cm]
\displaystyle \le 2^{\{s\}}\sum_{l\in\Z^n}\sum_{m=0}^{[s]} \begin{pmatrix}[s]\cr m\end{pmatrix}
\Big(\langle l\rangle^{m}|u_l|\langle k-l\rangle^{s-m}|v_{k-l}|\\[0.3cm]
 + 
\langle l\rangle^{m + \{s\}}|u_l|\langle k-l\rangle^{[s]-m}|v_{k-l}|\Big)
\\[0.3cm]
\displaystyle = 2^{\{s\}} \sum_{l\in\Z^n}\sum_{m=0}^{[s]} \begin{pmatrix}[s]\cr m\end{pmatrix}
\Big(\langle l\rangle^{m}|u_l|\langle k-l\rangle^{s-m}|v_{k-l}| + 
\langle l\rangle^{s-m}|u_l|\langle k-l\rangle^{m}|v_{k-l}|\Big).
\end{array}
\]
We have used the inequality $\left|\frac{a+b}{2}\right|^x\le \max\{a^x,b^x\}\le a^x+b^x$, where  $a,b\in\N$ and $x\ge 0$. Summing with respect to $k\in\Z^n$ 
we prove the claim.
\finishproof

A similar inequality can be obtained for the Sobolev $s$-norm of $uv$, but there appears an additional factor $2^{s/2}$ coming from the inequality $(a+b)^2\le 2(a^2+ b^2)$, which makes it useless for the estimates in Sect. \ref{sec:Gevrey}, because it changes the Gevrey constant at any step of the construction. 
 
To get rid of the sum in  Lemma \ref{lemma:product}, we  consider the modified norms
$$
P_s(u)\, =\, (s+1)^2S_s(u)\, ,\ s\ge 0\, ,\ u\in {\mathcal A}^s(\T^{n})\, .  
$$ 
If $f\in {\mathcal A}^{s+\tau}(\T^{n})$ and $\langle f\rangle =0$, and if $u\in {\mathcal A}^{s}(\T^{n})$ is a solution of the homological equation \eqref{eq:homological1} such that $\langle u\rangle =0$, then by Lemma \ref{lemma:homological} we obtain
\begin{equation}
\label{eq:homological2}
P_{s}(u)\,=(s+1)^2S_{s}(u)\, \le\,  \frac{(s+\tau+1)^2}{\kappa}\, S_{s+\tau}(f) =  \frac{1}{\kappa}\, P_{s+\tau}(f)\, . 
\end{equation}
Moreover, for any $u,v\in {\mathcal A}^s(\T^{n})$ we obtain from Lemma \ref{lemma:product} the following estimate
\begin{equation}
\label{eq:product}
\begin{array}{lcrr}
\displaystyle P_s(uv)\, \le\, 2 \sum_{m=0}^{[s]} \frac{(s+1)^2}{(m+1)^2(s-m+1)^2}\, \\[0.5cm]
\displaystyle \times\, \begin{pmatrix}[s]\cr m\end{pmatrix}
\big[P_{s-m}(u)P_m(v) + P_{s-m}(v)P_m(u)\big]\\[0.5cm]
\displaystyle  \le\,  \widetilde C\, \sup_{0\le m\le [s]}\, \left\{\begin{pmatrix}[s]\cr m\end{pmatrix}
\big[P_{s-m}(u)P_m(v) + P_{s-m}(v)P_m(u)\big]\right\}
\, ,
\end{array}
\end{equation}
where $\widetilde C=16\sum_{q=1}^{\infty} q^{-2}= 8\pi^2/3$. Another usefull property of the norm $P_s(\cdot)$, $s\ge 0$,  is that 
\[
P_s(\partial^\alpha u) \le P_{s+|\alpha|}(u) 
\]
for any $\alpha\in\N^n$ and $u\in C^\infty(\T^n)$. 

For any $p\in\N$ and $u\in C^\infty(\T^n)$ we set
\[
\displaystyle Q_p(u):= \sup_{ |\alpha|=p}\, \sup_{\theta\in\T^n}|\partial_\theta^\alpha u(\theta)|
\]
\begin{lemma}\label{lemma:gevrey-estimates}
There is a positive constant  $C_0=C_0(n)$ depending only on the dimension $n$ such that  
\begin{equation}
\label{eq:gevrey-estimates}
Q_{[s]}(u) \le P_s( u) \, 
\le\,  C_0\,  (2en)^{[s]}\, \left( Q_{[s]+n+2}(u) + Q_0(u)\right)
\end{equation}
for any $u\in C^\infty(\T^n)$ and $s\ge 0$. 
\end{lemma}
{\em Proof.} We have 
\[
\displaystyle P_s( u) \le C'_0(n)(1+s)^2\, \sup_{k\in \Z^n}\big((1+|k|)^{[s]+n+2}|u_k|\big)
\]
where $C^\prime_0(n):= \sum_{k\in\Z^n}(1+|k|)^{-n-1}$. 
Integrating by parts we get for any $p\in\N$ and any $k\neq 0$ the inequality
\[
\begin{array}{lcrr}
\displaystyle (1+|k|)^{p}|u_k|\, \le\,  (2n)^p\, \sup_{1\le j\le n}\, (|k_j|^p|u_k|)\, \le \, 
(2n)^{p}\, \sup_{|\alpha|=p}\, \sup_{\theta\in\T^n}\, |\partial_\theta^\alpha u(\theta)| \, .
\end{array}
\]
Moreover, $(1+s)^2\le 2e^{1+s}$, and 
we obtain  the second inequality in \eqref{eq:gevrey-estimates} with $C_0=2e^2(2n)^{n+2}C^\prime_0$. The proof of the first one is straightforward.
\finishproof 

Consider now the functions $b_\alpha$ given by \eqref{eq:hamiltonian}. 
\begin{lemma}\label{lemma:gevrey-estimates1}
We have 
\begin{equation}
\label{eq:gevrey-coefficients2}
P_s( b_\alpha) \le  \tilde L_0 L_1^s\,    L_2^{|\alpha|-1}\,  \Gamma(\rho s+(\rho-1)(|\alpha|-2)+1) 
\end{equation}
for any $s\ge 0$ and any $\alpha\in \N^n$ with a length  $|\alpha|\ge 2$, 
where the Gevrey constants $ L_1\ge 1$ and $ L_2\ge 1$  are equivalent to the corresponding Gevrey constants in 
\eqref{eq:gevrey-coefficients} and $\tilde L_0$ is equivalent to $L_0L_1^{n+2}$. 
\end{lemma}
Recall that the positive constant $\tilde L$ is equivalent to $L$ if there is $c(n,\rho,\tau)>0$ such that 
$\tilde L=c(n,\rho,\tau)L$. \\

{\em Proof.} 
Using  Lemma \ref{lemma:gevrey-estimates} and  \eqref{eq:gevrey-coefficients}, we get
\[
P_s( b_\alpha) \le   L_0 L_1^{s+n+2}  L_2^{|\alpha|-1}
\Gamma(\rho ([s]+n+2)+1) \Gamma((\rho-1)|\alpha|+1), 
\] 
where $L_0$, $L_1$ and $L_2$ are equivalent to the corresponding constants in \eqref{eq:gevrey-coefficients}. 
Note that the function $\Gamma(t)$ is increasing in the interval $[3/2, +\infty)$ and that 
$x^p\, \le\,  e^x\,  p\, !$ for any  $x\ge 0$ and $p\in \N$. Then 
using \eqref{eq:gamma-beta}, we obtain 
\[
\begin{array}{lcrr}
P_s( b_\alpha) \le   L_0 L_1^{s+n+2}\,  L_2^{|\alpha|-1} 
\Gamma(\rho s+(\rho-1)|\alpha|+\rho(n+2)+2)\\[0.5cm]
\le 
e^{\rho(n+4)}p!\, L_0 L_1^{n+2} (e^\rho L_1)^s\,  (e^{\rho-1} L_2)^{|\alpha|-1}\,  
\Gamma(\rho s+(\rho-1)(|\alpha|-2)+1)\, ,
\end{array}
\]
where $p=([\rho]+1)(n+4)$. This implies \eqref{eq:gevrey-coefficients2}. \finishproof

\section{Deriving the Homological Equation}\label{sec:homological} 
We turn now to the construction of the function $g$. The idea is to write explicitly the corresponding Taylor series and  to  prove certain Gevrey estimates for them  and then to use a Borel  extension theorem in Gevrey classes. 
Let us expand $g$ in Taylor series with respect to $I$ at $I=0$,
\begin{equation}
\label{eq:g}
g(\theta,I)\,  \sim\,  \sum_{m=2}^{\infty}\, g_m(\theta,I)\, ,\quad g_m(\theta,I) = \sum_{|\alpha|=m}g_{m,\alpha}(\theta)I^\alpha\, .
\end{equation}
Then we have formally
\[
\begin{array}{clr}
H(\theta, I + \partial g/\partial \theta (\theta, I)) \\[0.3cm]
\displaystyle{= \langle \omega, I\rangle + \sum_{m=2}^\infty 
\langle\omega,\frac{\partial g_m}{\partial \theta}(\theta, I) \rangle
+ \sum_{|\alpha|\ge 2}b_\alpha (\theta)\left(I + \sum_{k=2}^\infty \frac{\partial g_k}{\partial \theta}(\theta, I) \right)^\alpha}.
\end{array}
\]
We use the the following notations. For any $a=(a_1,\ldots,a_n)\in\C^n$ and $\alpha = (\alpha_1,\ldots,\alpha_n)\in\N^n$, we denote by $a^\alpha$ the product 
$a^\alpha := a_1^{\alpha_1}\cdots a_n^{\alpha_n}$, where by convention $z^0=1$ for any $z\in\C$. 
Let $$a_k=\{(a_{k,1},\ldots,a_{k,n})\in \C^n:\, k\in \N\}, $$ 
be a sequence in $\C^n$. Fix $\alpha\in \N^n$ of length   $|\alpha|\ge 2$, and 
recall  the following power series expansion in $\C[[X]]$ (c.f. (4.7) in \cite{G-P})
\[
\begin{array}{clr}
\displaystyle\left(\sum_{k=1}^\infty a_k X^k\right)^\alpha := 
\left(\sum_{k=1}^\infty a_{k,1} X^k\right)^{\alpha_1}
\cdots \left(\sum_{k=1}^\infty a_{k,n} X^k\right)^{\alpha_n}  \\[0.5cm]
\displaystyle= \sum_{m=|\alpha|}^\infty A_{\alpha,m} X^m\, , 
\end{array}
\]
where 
\[
\displaystyle
A_{\alpha,m} = \sum \, \frac{\alpha !}{\alpha^1 !\cdots \alpha^{m-1} !} \, 
a_1^{\alpha^1}\cdots a_{m-1}^{\alpha^{m-1}}\, ,
\]
and the sum is over the set $\N(\alpha,m)$  of all  multi-indices 
\[
(\alpha^1,\ldots \alpha^{m-1})\in \underbrace{\N^n\times\cdots\times\N^n}_{m-1}
\]
such that
\begin{equation}
\label{eq:index-set}
\left\{
\begin{array}{clr}
\alpha^1 +\cdots + \alpha^{m-1} = \alpha \, ,\\[0.3cm]
1\cdot|\alpha^1| + 2\cdot|\alpha^2| +\cdots + (m-1)\cdot|\alpha^{m-1}| = m .
\end{array}
\right.
\end{equation}
Notice that if $\alpha^1 +\cdots + \alpha^j = \alpha $, $1\cdot|\alpha^1| +\cdots + j\cdot|\alpha^j| = m$,  and $\alpha^j\neq 0$, then $j\le m-1$   since $|\alpha|\ge 2$. 
Hence, for any $\alpha\in \N^n$ with length $|\alpha|\ge 2$ we obtain 
\[
\left(I +\sum_{k=2}^\infty \frac{\partial g_k}{\partial \theta}(\theta, I) \right)^\alpha  
= \sum_{m=|\alpha|}^\infty A_{\alpha,m}(\theta,I)\, ,
\]
where $A_{\alpha,m}(\theta,I)$ is a homogeneous polynomial with respect to $I$ of degree $m$ of the form
\[
\begin{array}{lcrr}
\displaystyle
A_{\alpha,m}(\theta,m)\,  =\,  \sum \, \frac{\alpha !}{\alpha^1 !\alpha^2 !\cdots \alpha^{m-1} !} \, 
I^{\alpha^1}\\[0.5cm]
\displaystyle \times\, \left(\frac{\partial g_2}{\partial \theta}(\theta, I)\right)^{\alpha^2}\cdots \left(\frac{\partial g_{m-1}}{\partial \theta}(\theta, I)\right)^{\alpha^{m-1}}\, ,
\end{array}
\]
and the sum is taken over the set of multi-indices $\N(\alpha,m)$. 
Summing with respect to $\alpha$ we get formally
\begin{equation}
\label{eq:new-hamiltonian}
\begin{array}{lcrr}
\displaystyle H(\theta, I + \partial g/\partial \theta (\theta, I))  \\[0.5cm]
\displaystyle = \, \langle \omega, I\rangle\,   +\,  \sum_{m=2}^\infty \left(
\langle\omega,\frac{\partial g_m}{\partial \theta}(\theta, I) \rangle
\, +\,  B_m(\theta,I)\right)\, ,
\end{array}
\end{equation}
where 
$B_m(\cdot, I)$  is a  homogeneous 
 polynomial of degree $m\ge 2$ with respect to $I$  of the form 
\begin{equation}
\begin{array}{lcrr}
\displaystyle
B_{m}(\theta,I)\ = \  \sum \, \frac{\alpha !}{\alpha^1 !\alpha^2 !\cdots \alpha^{m-1} !} \,  b_\alpha(\theta)\,
I^{\alpha^1} \\[0.5cm]
\displaystyle\times\, \left(\frac{\partial g_2}{\partial \theta}(\theta, I)\right)^{\alpha^2}\cdots \left(\frac{\partial g_{m-1}}{\partial \theta}(\theta, I)\right)^{\alpha^{m-1}}\, .
\end{array}
\label{eq:remainder}
\end{equation}
The index set of the sum above is 
\[
\N(m):= \bigcup_{2\le |\alpha|\le m} \N(\alpha,m)
\]
and it consists of all the multi-indices 
\[(\alpha^1,\alpha^2, \ldots \alpha^{m-1})\in \underbrace{\N^n\times\cdots\times\N^n}_{m-1}\]
such that 
\[
1\cdot|\alpha^1| +2\cdot|\alpha^1| +\cdots + (m-1)\cdot|\alpha^{m-1}| = m\, . 
\]
Note    that  $\alpha = \alpha^1 +\alpha^2+ \cdots + \alpha^{m-1}$ in (\ref{eq:remainder}).

For any $m\ge 2$ we  obtain from \eqref{eq:normal-form} the following homological equation 
\begin{equation}
\label{eq:homological}
{\mathcal L}_\omega g_m(\theta,I) + B_m(\theta,I) = R_m(I)\, ,\quad \langle g_m(\cdot,I)\rangle = 0\, ,
\end{equation}
where
\[
R_m(I) = \langle  B_m(\cdot,I)\rangle = \frac{1}{(2\pi)^n}\int_{\T^n} B_m(\theta,I)\, d\theta\, .
\]
Denote by $(B_m)_k(I)$, $k\in \Z^n$, the Fourier coefficients of $B_m(\theta,I)$. Then
$(B_m)_0(I) = R_m(I)$ is a  homogeneous 
 polynomial of degree $m$ with respect to $I$. Since $\omega$ satisfies (\ref{eq:sdc})
\begin{equation}
\label{eq:solution}
g_m(\theta,I) = - \sum_{k\in\Z^n\setminus \{0\}} \frac{(B_m)_k(I)}{i\langle\omega,k\rangle}
e^{i\langle k,\theta\rangle}\, ,
\end{equation}
and it is a  homogeneous 
 polynomial of degree $m$ with respect to $I$ with smooth coefficients $g_{m,\alpha}(\theta)$, $|\alpha|=m$. Our aim is to obtain Gevrey estimates for $g_{m,\alpha}(\theta)$.

\section{Gevrey estimates}\label{sec:Gevrey}
We are going to show that there are positive constants $C_1$ and $C_2$ depending on the constants $\tilde L_0$,  $L_1$ and $L_2 $ in \eqref{eq:gevrey-coefficients2} 
such that for any $\alpha\in\N^n$ with length  $m=|\alpha|\ge 2$ and for any  $\beta\in\N^n$ we have 
\begin{equation}\label{eq:gevrey-estimates-coeff}
\displaystyle \sup_{\theta\in\T^n}\, |\partial^\beta g_{m,\alpha}(\theta)|\, \le \, C_1^{|\beta|} C_2^{|\alpha|-1} \beta !^\rho \alpha ! ^{\mu -1},
\end{equation}
where  $\mu =\rho(\tau+1)+1$ (see the statement of Theorem \ref{theo:gevrey}).
Consider for any $m\ge 2$ the solution ($g_m$, $R_m$) of the homological equation (\ref{eq:homological}), where $\omega$ is $(\kappa,\tau)$-Diophantine, $0<\kappa\le 1$ and $\tau>n-1$. 
Denote the unit poly-disc in $\C^n$ by $\D^n$, i.e.  $I=(I_1,\ldots,I_n)\in\C^n$ belongs to $\D^n$ if  $|I_j|\le 1$ for any $1\le j \le n$. 
\begin{prop}
\label{prop:main-estimates} 
There is   $A_0=A_0(n, \rho,\tau)\ge 1$ depending only  on  $n$, $\rho$ and $\tau$, 
such that  
\begin{equation}
\mbox{for}\ C_1 =  e^{\rho}L_1\ \mbox{and for any}\ C_2\ge  \frac{1}{\kappa}\, A_0 \tilde L_0 L_1^{\tau+2} L_2 
\label{eq:theconstant}
\end{equation}
the following estimates hold  
\begin{equation}
\label{eq:gevrey-estimates1}
\sup_{I\in \D^n}P_s(B_m(\cdot,I))\le  B_0 \tilde L_0 L_1^2 L_2 C_1^{s}C_2^{m -2}
\Gamma\big(\rho s +(\mu-1)(m-2)\big)
\end{equation}
for $m\ge 3$ and  any $s\in \R_+$, and 
\begin{equation}
\label{eq:gevrey-estimates2}
\sup_{I\in \D^n} P_s(g_m(\cdot,I)) \le  C_1^{s}C_2^{m -1}
\Gamma\big(\rho s + (\mu-1)(m-1)-\rho \big),
\end{equation}
for $m\ge 2$ and  any $s\in \R_+$, where $\tilde L_0$, $L_1$ and $L_2$ are the corresponding Gevrey constants in  \eqref{eq:gevrey-coefficients2} and $B_0=B_0(n,\rho,\tau)\ge 1$. 
\end{prop}
Note that $(\mu-1)(m-1)-\rho \ge  \rho(\tau+1)-\rho>1$.
We are going to prove Proposition \ref{prop:main-estimates} by recurrence with respect to $m\ge 2$.  
For $m=2$ we obtain from (\ref{eq:homological}) the equation 
\[
{\mathcal L}_\omega g_2(\theta,I) = R_2(I)-B_2(\theta,I)\, . 
\]
Moreover, \eqref{eq:index-set} and \eqref{eq:remainder} imply
\[
B_2(\theta,I)= \sum_{|\alpha|=2}\, b_\alpha(\theta) I^\alpha \, .
\]
For any $I\in\D^n$ we have by \eqref{eq:gevrey-coefficients2}
\[
\begin{array}{lcrr}
\displaystyle  P_s(B_2(\cdot,I)) \le \sum_{|\alpha|=2}P_s(b_\alpha) 
\le  n^2 \,  \tilde L_0 L_1^s L_2 \, \Gamma(\rho s + 1). 
\end{array}
\] 
Then using  \eqref{eq:homological2} we obtain 
\[
\begin{array}{lcrr}
\displaystyle P_s(g_2(\cdot,I)) \le \frac{1}{\kappa}P_{s+\tau}(B_2(\cdot,I))
\le n^2 \,  L_1^s  \left(\frac{1}{\kappa}\tilde L_0L_1^{\tau} L_2\right)\, \Gamma(\rho s + \rho\tau +1).\\[0.3cm] 
\end{array}
\]
On the other hand, 
\[
\Gamma(\rho s + \rho\tau +1) = (\rho s + \rho\tau )\Gamma(\rho s + \rho\tau) 
\le e^{\rho s + \rho\tau}\Gamma(\rho s + \rho\tau)\, ,
\]
and we obtain
\[
\displaystyle P_s(g_2(\cdot,I)) \le  C_1^s  C_2 \Gamma(\rho s + (\mu-1)-\rho )
\]
for any $I\in \D^n$, where $C_1 = e^\rho L_1$, $C_2\ge \frac{1}{\kappa}A_0\tilde L_0L_1^{\tau}L_2$ and $A_0\ge n^2e^{\rho\tau}$. 

Fix $m\ge 3$ and suppose that the estimates (\ref{eq:gevrey-estimates2}) hold for any $p<m$ and any $s\ge 0$.
We are going first to estimate $P_s(B_m(\cdot,I))$, $I\in\D^n$, for any $s\ge 0$. Using the inductive assumption we get 
\begin{lemma}\label{lemma:estimate1}
Let $p\ge 1$ and $2\le m_k\le m-1$, where $k=1,\ldots,  p$. 
Set $$M_p=m_1+\cdots+m_p-p.$$
 Then for any $\delta\in (0,\mu-1)$ there is a constant $ C_0(\delta,\mu)\ge 1$ such that
\begin{equation}\label{eq:product-estimates1}
\begin{array}{lcrr}
\displaystyle
\sup_{I\in\D^n}P_s\left( \frac{\partial g_{m_1}}{\partial \theta_{j_1}} (\cdot,I)\ldots \frac{\partial g_{m_p}}{\partial \theta_{j_p}}(\cdot,I)\right) \le C_0(\delta,\mu)^{p-1} C_1^{p+s} C_2^{M_p} \\ [0.5cm]
\displaystyle \times\, \left(\frac{M_p!}{(m_1-1)!\cdots (m_p-1)!}\right)^{1+\delta-\mu}\, 
\Gamma\big(\rho s + (\mu-1)M_p \big).
\end{array}
\end{equation}
\end{lemma}
\begin{rem}
\label{re:M}
If $p\ge 2$ then
\[
\frac{M_p!}{(m_1-1)!\cdots (m_p-1)!}\ge M_p\, . 
\]
\end{rem}
\begin{rem}
\label{re:delta}
Recall that  $\mu =\rho(\tau +1)+1 > 3$. We shall fix later $\delta=\mu-2>1$. 
\end{rem}
\noindent
{\em Proof of Lemma \ref{lemma:estimate1}}. We are going to prove \eqref{eq:product-estimates1} by induction with respect to $p\ge 1$. 
For $p=1$, we have
\begin{equation}
\label{eq:derivative}
P_s\left( \frac{\partial g_{m_1}}{\partial \theta_{j}}\right) \le P_{s+1}(g_{m_1}) \le  C_1^{s+1}C_2^{m_1 -1}
\Gamma\big(\rho s + (\mu-1)(m_1-1) \big)
\end{equation}
in view of  \eqref{eq:gevrey-estimates2}. 
Set 
\[
\displaystyle F_p(\theta,I)=\frac{\partial g_{m_1}}{\partial \theta_{j_1}} (\theta,I)\ldots \frac{\partial g_{m_p}}{\partial \theta_{j_p}}(\theta,I). 
\]
Now take   $p=2$ and   $2\le m_1,\, m_2\le m-1$, and fix  $I\in \D^n$. Using \eqref{eq:product} and \eqref{eq:derivative} we obtain
\[
\begin{array}{lcrr}
\displaystyle
P_s\left( F_2(\cdot,I)\right)
 \le\,  \widetilde C\, C_1^{s+2} C_2^{m_1+m_2-2}\,  \max _{0\le q\le [s]}\Big[
\begin{pmatrix}[s]\cr q\end{pmatrix} \\[0.5cm]
 \times\, \Big(\Gamma(\rho(s-q) + (\mu-1)(m_1-1)  )\Gamma(\rho q + (\mu-1)(m_2-1)  )\\[0.5cm]
 + \,  \Gamma(\rho(s-q) + (\mu-1)(m_2-1)  )\Gamma(\rho q+ (\mu-1)(m_1-1)  )\Big)\Big]\, .   
 \end{array}
\]
On the other hand, 
\[
\begin{array}{lcrr}
\Gamma(\rho(s-q) + (\mu-1)(m_1-1)  )\Gamma(\rho q + (\mu-1)(m_2-1)  )\\[0.5cm]
= \Gamma(\rho s + (\mu-1)(m_1+m_2-2) ) \\[0.5cm]
\times B(\rho(s-q) + (\mu-1)(m_1-1) ,  \rho q + (\mu-1)(m_2-1)  ),  
\end{array}
\]
where $B(x,y)$, $x,y>0$, is the  Beta function  (\ref{eq:BetaFunction}). Recall that $B(x,y)$ is decreasing with respect to both variables  $x>0$ and $y>0$. Then using  \eqref{eq:cauchy} we get for any $\delta\in (0,\mu-1)$ the inequalities
\[
\begin{array}{lcrr}
B(\rho(s-q) + (\mu-1)(m_1-1) ,  \rho q + (\mu-1)(m_2-1)  )\\[0.5cm]
\le 
B(2\rho(s-q) + \delta,  2\rho q + \delta  )^{1/2}\\[0.5cm]
\times\, B(2(\mu-1)(m_1-1) -\delta,  2(\mu-1)(m_2-1) -\delta )^{1/2}\\[0.5cm]
\le B(2(s-q) + \delta,  2q + \delta  )^{1/2}\\[0.5cm]
\times\, 
B(2(\mu-1-\delta)(m_1-1) +\delta,  2(\mu-1-\delta)(m_2-1) + \delta )^{1/2}.
\end{array}
\]
Moreover, Lemma \ref{lemma:gamma} implies 
\begin{equation}\label{eq:estimate-of-B}
\begin{array}{lcrr}
B(2(s-q) + \delta,  2q + \delta  )  \\[0.5cm]
 \displaystyle  \le B(2([s]-q) + \delta,  2  q + \delta  )
\le C^\prime (\delta) \begin{pmatrix}[s]\cr q\end{pmatrix}^{-2}.
 \end{array}
 \end{equation}
as well as 
\begin{equation}\label{eq:estimate-of-B2}
\begin{array}{lcrr}
B(2(\mu-1-\delta)(m_1-1) +\delta,  2(\mu-1-\delta)(m_2-1) + \delta )\\[0.5cm]
\le 
\displaystyle C^\prime (\delta,\mu) \begin{pmatrix}m_1+m_2-2\cr m_1-1\end{pmatrix}^{2(1+\delta-\mu)}.
 \end{array}
 \end{equation}
Hence, for any non-negative integer $ 0\le q\le [s]$ we have 
\[
\begin{array}{lcrr}
\displaystyle \begin{pmatrix}[s]\cr q\end{pmatrix} B(\rho(s-q) + (\mu-1)(m_1-1) ,  \rho q + (\mu-1)(m_2-1)  )\\[0.5cm]
\displaystyle
\le C''(\delta,\mu)\,  \begin{pmatrix}m_1+m_2-2\cr m_1-1\end{pmatrix}^{1+\delta-\mu}\, .
 \end{array}
\]
In the same way we estimate the quantity
\[
\Gamma(\rho(s-q) + (\mu-1)(m_2-1)  )\Gamma(\rho q + (\mu-1)(m_1-1)  ). 
\]
Finally we obtain  
\[
\begin{array}{lcrr}
\displaystyle P_s\left( F_2(\cdot,I)\right)\,  \le \, 
 C_0(\delta,\mu)\, C_1^{s+2} C_2^{m_1+m_2-2}\\[0.5cm]
 \displaystyle \times\, \Gamma(\rho s + (\mu-1)(m_1+m_2-2) )\left(\frac{(m_1+m_2-2)!}{(m_1-1)!(m_2-1)!}\right)^{1+\delta-\mu}\, ,
 \end{array}
\]
where
\begin{equation}
\label{eq:crho1}
C_0(\delta, \mu)\, =\, \max\{ 2\, \widetilde C\, C''(\delta,\mu), 1\}\ge 1. 
\end{equation}
The proof follows by recurrence with respect to $p$  setting $F_{p}=F_{p-1}\frac{\partial g_{m_p}}{\partial \theta_{j_p}}$ and then using  \eqref{eq:product-estimates1} for $F_{p-1}$, \eqref{eq:derivative}, and the above argument. 
At any step the constant $C_0(\delta, \mu)$ is given by (\ref{eq:crho1}). 
\finishproof\\
In the same way as above, using (\ref{eq:gevrey-coefficients2}), we get
\begin{lemma}\label{lemma:estimate2}
Let  $p\ge 1$ and $2\le m_k\le m-1$, where $k=1,\ldots,  p$,  and let $\alpha\in \N^n$ with  $|\alpha|\ge 2$. 
Set $M_p=m_1+\cdots+m_p-p$ and  $C_1= e^\rho L_1$.  Then for any $0<\delta < \mu-1$ and   $I\in\D^n$  we have
\[
\begin{array}{lcrr}
\displaystyle
P_s\left(b_\alpha \frac{\partial g_{m_1}}{\partial \theta_{j_1}} (\cdot,I)\ldots \frac{\partial g_{m_p}}{\partial \theta_{j_p}}(\cdot,I)\right)\\[0.5cm] 
\displaystyle \le \
K_0  \tilde L_0 (e^{\mu-1}L_2)^{|\alpha|-1} C_0^{p-1}  C_1^{p+s}  C_2^{M_p}\, \Gamma\big(\rho s + (\mu-1)(M_p+|\alpha|-2)  \big) 	 \\[0.5cm] 
\displaystyle \times \left(\frac{M_p!}{(m_1-1)!\cdots (m_p-1)!}\right)^{1+\delta-\mu}\, \begin{pmatrix} M_p +|\alpha|-2 \cr |\alpha|-2 \end{pmatrix}^{1+\delta-\mu}\,  , 
\end{array}
\]
where  $C_0=C_0(\delta,\mu)\ge 1$ is given by \eqref{eq:crho1}  and $K_0 =K_0(n,\delta,\mu)\ge 1$.
\end{lemma}
\noindent
{\em Proof}. 
To simplify the notations we will write below $M$ instead of $M_p$. Set as above $\displaystyle F_p(\theta,I)=\frac{\partial g_{m_1}}{\partial \theta_{j_1}} (\theta,I)\ldots \frac{\partial g_{m_p}}{\partial \theta_{j_p}}(\theta,I)$ and fix $I\in\D^n$. First suppose that  $|\alpha|=2$. Using \eqref{eq:product} and  Lemma \ref{lemma:estimate1}, we obtain 
\[
\begin{array}{lcrr}
\displaystyle
P_s\left(b_\alpha \, F_p(\cdot,I)\right)\, 
 \, \\[0.5cm] 
\displaystyle \le\,   \widetilde C  C_0^{p-1}C_1^{p}  C_2^{M}
\,  \left(\frac{M_p!}{(m_1-1)!\cdots (m_p-1)!}\right)^{1+\delta-\mu} \\[0.5cm] 
\displaystyle
\times\,    \max_{0\le q\le [s]}\Big[
\begin{pmatrix}[s]\cr q\end{pmatrix}   \Big(  P_q(b_\alpha)\Gamma\big(\rho(s-q) + (\mu-1) M  \big)C_1^{s-q}\\[0.5cm] 
+\ P_{s-q}(b_\alpha) \Gamma(\rho q + (\mu-1) M  \big)C_1^{q}\Big)\Big]\, .  
 \end{array}
\]
For $q=0$ we have 
\[
P_0(b_\alpha)\Gamma\big(\rho s + (\mu-1) M  \big)\le  \tilde L_0 L_2\Gamma\big(\rho s + (\mu-1) M  \big).
\]
For $q\ge 1$ the estimates (\ref{eq:gevrey-coefficients2}) imply
\[
P_q(b_\alpha)\le \tilde L_0L_1^qL_2\Gamma\big(\rho q + 1)\le \tilde L_0(e^\rho L_1)^qL_2\Gamma\big(\rho q  \big)= \tilde L_0 C_1^qL_2\Gamma\big(\rho q  \big),
\]
and we obtain
\[
\begin{array}{lcrr}
P_q(b_\alpha)\Gamma\big(\rho(s-q) + (\mu-1) M  \big) \\[0.3cm] 
\le 
\tilde L_0 C_1^qL_2\Gamma\big(\rho s + (\mu-1) M  \big)B(\rho q, \rho(s-q) + (\mu-1) M ).
\end{array}
\]
On the other hand, $\rho\ge 1$, $(\mu-1) M \ge \mu-1=\rho(\tau+1)>2$, and  $B(x,y)$ is decreasing with respect to both variables $x>0$ and $y>0$, hence,  we get
\[
\begin{array}{lcrr}
B(\rho q, \rho(s-q) + (\mu-1) M )\\[0.3cm] 
\displaystyle \le B(q, [s]-q + 1 )=\frac{(q-1)!([s]-q )!}{[s]!} < \begin{pmatrix}[s]\cr q\end{pmatrix}^{-1}.
\end{array}
\]
In the same way we estimate the second term of the sum above. For $0\le q<[s]$ we use the same argument as above. For $q=[s]$ we get 
$$P_{\{s\}}(b_\alpha)\le \tilde L_0 L_1^{\{s\}}L_2\Gamma(\rho{\{s\}}+1)\le \tilde L_0 L_1^{\{s\}}L_2\Gamma(\rho+1),$$
 since $\rho+1\ge 2$ (see the argument below), 
and we obtain
\[
\begin{array}{lcrr}
P_{\{s\}}(b_\alpha)\Gamma\big(\rho [s] + (\mu-1) M  \big) \\[0.3cm] 
\le 
\Gamma(\rho+1)\tilde L_0 L_1^{\{s\}} L_2\Gamma\big(\rho s + (\mu-1) M  \big).
\end{array}
\]
This proves the claim for $|\alpha|=2$. 

Let $|\alpha|\ge 3$. Recall that $\Gamma(t)$ is increasing in the interval $[2,+\infty)$, 
$ \Gamma(t)\le 1$ for $t\in [1,2]$, and $\Gamma(1)=\Gamma(2)=1$ (see Sect. \ref{sec:Gamma}).
Hence, $ \Gamma(t_1)\le  \Gamma(t_2)$ if $1\le t_1\le t_2$ and $t_2\ge 2$.  Since  $\mu-1=\rho(\tau+1)>2$ and $|\alpha|-2\ge 1$, 
this allows us to replace $\rho-1$ by $\mu-1$ in (\ref{eq:gevrey-coefficients2}), and we obtain
\begin{equation}
\label{eq:gevrey-coefficients3}
\begin{array}{lcrr}
P_s(b_\alpha) \le \tilde L_0 L_1^s L_2^{|\alpha|-1}\Gamma(\rho s + (\mu-1)(|\alpha|-2)+1)\\[0.3cm]
\le \tilde L_0 (e^{\rho}L_1)^s (e^{\mu-1}L_2)^{|\alpha|-1}\Gamma(\rho s + (\mu-1)(|\alpha|-2)) 
\end{array}
\end{equation}
for any $\rho\ge 1$ and $s\ge 0$. 
Using  \eqref{eq:product} and  Lemma \ref{lemma:estimate1}, we obtain 
\[
\begin{array}{lcrr}
\displaystyle
P_s\left(b_\alpha \, F_p(\cdot,I)\right)\,  \\[0.5cm] 
\displaystyle   \le\,   \widetilde C \tilde L_0 C_0^{p-1}C_1^{s+p}(e^{\mu-1}L_2)^{|\alpha|-1} C_2^{M} \,
\left(\frac{M!}{(m_1-1)!\cdots (m_p-1)!}\right)^{1+\delta-\mu}\,  \,\\[0.5cm] 
\displaystyle
\times\, \max_{0\le q\le [s]}\Big[
\begin{pmatrix}[s]\cr q\end{pmatrix}      \big(  \Gamma\big(\rho q+(\mu-1)(|\alpha|-2) \big)\Gamma\big(\rho(s-q) + (\mu-1) M  \big)\\[0.5cm]
+\, \Gamma\big(\rho(s-q)  +(\mu-1)(|\alpha|-2))\Gamma(\rho q + (\mu-1) M  \big)\big)\Big]\, .  
 \end{array}
\]
Recall that $|\alpha|\ge 3$, hence $\rho q+(\mu-1)(|\alpha|-2)\ge \mu-1=\rho(\tau+1) >2$. We have
\[
\begin{array}{lcrr}
\displaystyle
\Gamma\big(\rho q  +(\mu-1)(|\alpha|-2)\big)\Gamma\big( \rho (s-q) + (\mu-1) M  \big)\\[0.5cm]
\displaystyle
=\, \Gamma\big(\rho s  +(\mu-1)(M+|\alpha|-2) \big)\\[0.5cm]
 \displaystyle \times\, 
B\big(\rho q  +(\mu-1)(|\alpha|-2), \rho (s-q) + (\mu-1) M  \big)\, .
\end{array}
\]
Using  \eqref{eq:cauchy} we get 
\[
\begin{array}{lcrr}
\displaystyle
B(\rho q  +(\mu-1)(|\alpha|-2), \rho (s-q) + (\mu-1) M  ) \\[0.5cm]
\le  B(2\rho q  +\delta, 2\rho (s-q) + \delta  )^{1/2}  \\[0.5cm]
\times B(2(\mu-1)(|\alpha|-2)-\delta,  2(\mu-1)M -\delta )^{1/2}. 
 \end{array}
\]
Moreover, using  Lemma \ref{lemma:gamma} we get as above 
\[
\begin{array}{lcrr}
\displaystyle
B\big(2\rho q  +\delta, 2\rho (s-q) + \delta \big)
 \le\,   B\big(2q  + \delta,  2([s]-q) + \delta \big) \\[0.5cm]
 \le\, C'(\delta) \begin{pmatrix}[s]\cr q\end{pmatrix}^{-2}\, ,
 \end{array}
\]
and 
\[
\begin{array}{lcrr}
\displaystyle
B(2(\mu-1)(|\alpha|-2)-\delta,  2(\mu-1)M -\delta )\\[0.5cm]
\le B(2(\mu-1-\delta)(|\alpha|-2)+\delta,  2(\mu-1-\delta)M+\delta )\\[0.5cm]
\displaystyle\le C'(\delta,\mu)\, 
\begin{pmatrix}M +|\alpha|-2\cr |\alpha|-2\end{pmatrix}^{2(1+\delta-\mu)}\, . 
 \end{array}
\]
In the same way we estimate the second term. 
This proves the Lemma taking  $K_0= K_0(n,\delta,\mu)\ge 1$ sufficiently large. 
\finishproof

From now on we fix $\delta =\mu-2= \rho(\tau+1)-1\ge \tau >1$. 
We return  to the proof of (\ref{eq:gevrey-estimates1}) and (\ref{eq:gevrey-estimates2}). 
First we shall estimate $P_s(B_{m}(\cdot,I))$ for $I\in \D^n$ and  $m\ge 3$. 
In view of  (\ref{eq:remainder}) we obtain 
\[
\displaystyle P_s(B_{m}(\cdot,I))\
   \le\   \sum_{2\le |\alpha|\le m} Q_\alpha (I),
\]
where
\begin{equation}
\label{eq:Q}
\begin{array}{lcrr}
\displaystyle Q_\alpha (I)
   =   \sum_{(\alpha^1,\ldots,\alpha^{m-1})\in \N(\alpha,m)} \ \frac{\alpha !}{\alpha^1 !\cdots \alpha^{m-1} !}\,  \\[0.5cm]
   \displaystyle \times\,  P_s\left( b_\alpha\,
\left(\frac{\partial g_2}{\partial \theta}(\cdot, I)\right)^{\alpha^2}\cdots \left(\frac{\partial g_{m-1}}{\partial \theta}(\cdot, I)\right)^{\alpha^{m-1}}\right).
\end{array}
\end{equation}
Consider more closely the index set $\N(\alpha,m)$, where $|\alpha|\ge 2$. Recall from \eqref{eq:index-set} that 
$(\alpha^1,\ldots, \alpha^{m-1})\in (\N^n)^{m-1}$ belongs to $\N(\alpha,m)$ if and only if 
\[
\left\{
\begin{array}{lcrr}
\alpha^1 +\cdots + \alpha^{m-1} = \alpha,  \ \mbox{and}\\ 
1\cdot|\alpha^1| + 2\cdot|\alpha^2| +\cdots + (m-1)\cdot|\alpha^{m-1}| = m .
\end{array}
\right.
\]
Set 
\[
\begin{array}{lcrr}
\N_0(\alpha,m):=\{(\alpha^1,\alpha^2,\ldots,\alpha^{m-1})\in \N(\alpha,m):\, \alpha_1=\alpha\}\, , \\[0.3cm]
\N_1(\alpha,m):=\{(\alpha^1,\alpha^2,\ldots,\alpha^{m-1})\in \N(\alpha,m):\, |\alpha_1-\alpha|=1\}\, ,\\[0.3cm]
\N^\ast(\alpha,m):=\{(\alpha^1,\alpha^2,\ldots,\alpha^{m-1})\in \N(\alpha,m):\, |\alpha_1-\alpha|\ge 2\}\, , 
\end{array}
\]
and denote the corresponding sums in \eqref{eq:Q} by $Q^0_\alpha (I)$, $Q^1_\alpha (I)$ and $Q^\ast_\alpha (I)$ respectively. \\

\noindent
{\em  1. Estimate of $Q^0_\alpha (I)$}. 
The set $\N_0(\alpha,m)$ contains only one element,  $|\alpha|=m\ge 3$,  and  by \eqref{eq:gevrey-coefficients3} we get
\[
\begin{array}{lcrr}
Q^0_\alpha (I) \le P_s(b_\alpha) \le \tilde L_0  (e^{\rho}L_1)^{s} (e^{\mu-1}L_2)^{m-1} 
 \Gamma(\rho s + (\mu-1)(m-2)  )\\[0.3cm]
\le\,   2^{1-m} C_1^{s}C_2^{m-1}  \Gamma(\rho s + (\mu-1)(m-2)  )
 \end{array}
\]
for $C_2\ge 2 e^{\mu-1}L_2$ and $C_1=e^\rho L_1$. \\

\noindent
{\em  2. Estimate of $Q^1_\alpha (I)$}. 
Notice that the cardinality of $\N_1$ is $\# \N_1(\alpha,m)\le n$. Indeed, if $(\alpha^1,\alpha^2,\ldots,\alpha^{m-1})\in \N_1(\alpha,m)$, then we have $|\alpha^1|= |\alpha|-1\ge 1$ and $|\alpha^2| +\cdots + |\alpha^{m-1}|=1$. 
Hence, $\alpha^k=0$ for any $k\neq 1, m-|\alpha|+1$, and $|\alpha^k|=1$ for $k=m-|\alpha|+1$, which implies   $\# \N_1(\alpha,m)\le n$. Moreover, $\alpha!/\alpha^1 ! \le |\alpha|$. 

Fix $C_1=e^\rho L_1$ and  $C_2\ge 2e^{\mu-1}L_2$. 
Using  Lemma \ref{lemma:estimate2} with $p=1$, $m_1= m-|\alpha|+1$ and $M_1=m-|\alpha|$,  we get  
\[
\begin{array}{lcrr}
Q^1_\alpha (I) \le\,  |\alpha| K_0'\tilde L_0 (e^{\mu-1}L_2)^{|\alpha|-1} C_1^{s+1} C_2^{m-|\alpha|} 
 \Gamma(\rho s + (\mu-1)(m-2)  ) \\[0.3cm]
 \le\,  |\alpha|    (K_0'\tilde L_0 L_1L_2)   (e^{\mu-1}L_2)^{|\alpha|-2} C_1^{s} C_2^{m-|\alpha|} 
 \Gamma(\rho s + (\mu-1)(m-2)  )\\[0.3cm]
 \le\,  |\alpha|   2^{- |\alpha|}(K_0'\tilde L_0 L_1L_2) C_1^{s} C_2^{m-2} 
 \Gamma(\rho s + (\mu-1)(m-2)  )
 \, ,
 \end{array}
\]
where  $K_0'=K_0'(n,\rho,\mu) $ stands for different constants  
depending only on $n$, $\rho$ and $\mu$.  \\

\noindent
{\em  3. Estimate of $Q^\ast_\alpha (I)$}. 
Let $(\alpha^1,\alpha^2,\ldots,\alpha^{m-1})\in \N^\ast(\alpha,m)$. Set as above
\[
F:= \left(\frac{\partial g_2}{\partial \theta}(\cdot, I)\right)^{\alpha^2}\cdots \left(\frac{\partial g_{m-1}}{\partial \theta}(\cdot, I)\right)^{\alpha^{m-1}}.
\]
Notice that  the corresponding $p$ in  Lemma \ref{lemma:estimate2} is 
$$p=|\alpha_2|+\cdots +|\alpha_{m-1}|\ge 2. $$
Moreover, $p\le |\alpha|$ and 
\[
\begin{array}{lcrr}
M_p:=1\cdot|\alpha_1|+2\cdot|\alpha_2|+\cdots + (m-1)\cdot|\alpha_{m-1}|\\[0.3cm]
-|\alpha_1 +\alpha_2+\cdots+ \alpha_{m-1}|
= m-|\alpha| \ge 2\, .
\end{array}
\]
It follows from   Lemma \ref{lemma:estimate2} and Remark \ref{re:M} that 
\[
\begin{array}{lcrr}
\displaystyle P_s\left( b_\alpha\, F\right)
  \le\   K_0  \tilde L_0  \big(C_0e^{\mu-1}L_2\big)^{|\alpha|-1}  C_1^{s+|\alpha|} C_2^{m-|\alpha|} \\[0.5cm]
\times\,  \Gamma(\rho s + (\mu-1)(m-2)  )\,  \begin{pmatrix} m-2 \cr |\alpha|-2 \end{pmatrix}^{-1}\, (m-|\alpha|)^{-1}
\end{array}
\]
for any $I\in \D^n$. We have 
\[
\displaystyle  \begin{pmatrix}m-1 \cr |\alpha|-1 \end{pmatrix}\begin{pmatrix}m-2 \cr |\alpha|-2 \end{pmatrix}^{-1}(m-|\alpha|)^{-1} = \frac{m-1}{(|\alpha|-1)(m-|\alpha|)} \le 2. 
\]
Then using Lemma \ref{lemma:binom} we estimate  $Q_\alpha^\ast (I)$  by
\[
\begin{array}{lcrr}
\displaystyle  Q_\alpha^\ast (I)\le K_0 \tilde L_0  \big(C_0e^{\mu-1} L_2\big)^{|\alpha|-1} 
 C_1^{s+|\alpha|}  C_2^{m-|\alpha|} 
\, \Gamma(\rho s + (\mu-1)(m-2)  )	\\[0.5cm]
\displaystyle \times \begin{pmatrix}m-1 \cr |\alpha|-1 \end{pmatrix}\begin{pmatrix}m-2 \cr |\alpha|-2 \end{pmatrix}^{-1}(m-|\alpha|)^{-1}
\\[0.5cm]
\displaystyle \le (K_0'\tilde L_0L_1^2L_2)  \big(C_0 e^\rho e^{\mu-1} L_1L_2\big)^{|\alpha|-2} C_1^{s} C_2^{m-|\alpha|} 
\, \Gamma(\rho s + (\mu-1)(m-2)  )\, . 
\end{array}
\]
Hereafter $K_0'=K_0'(n,\rho,\mu) $ stands for different constants  
depending only on $n$, $\rho$ and $\mu$.
For   $C_2\ge   2C_0(\delta,\mu) e^\rho e^{\mu-1} L_1L_2 $ we obtain 
\[
\displaystyle  Q_\alpha^\ast (I)\le  2^{-|\alpha|}(K_0'C_0\tilde L_0L_1^2L_2) C_1^{s}  C_2^{m-2} 
\, \Gamma(\rho s +(\mu-1) (m-2)  )\, . 
\]
Taking into account the cases {\em 1. -  3.}  we obtain 
\[
\displaystyle  Q_\alpha (I)\le |\alpha|2^{-|\alpha|}  (K_0'\tilde L_0L_1^2L_2) C_1^{s}   C_2^{m-2} 
\, \Gamma(\rho s + (m-2)(\mu-1)  )\,  
\]
for any $I\in \D^n$, where 
\begin{equation}
C_1=e^\rho L_1 \ \mbox{and}\ C_2\ge   2C_0(\delta,\mu)e^{\rho}e^{\mu-1} L_1^2L_2 .
\label{eq:constants1}
\end{equation}
Set 
\[
B_0:=  K_0' \sum_{p=0}^\infty (p+1)^{n+1}2^{-p}.
\] 
Then for any $s\ge 0$, $m\ge 3$, and   $I\in \D^n$ we obtain 
\[
\displaystyle P_s(B_{m}(\cdot,I))\
 \le\, B_0 \tilde L_0L_1^2L_2 C_1^{s} C_2^{m-2}\Gamma(\rho s + (\mu-1)(m-2)  )\, ,
\]
which proves \eqref{eq:gevrey-estimates1}.

This implies
\[
|R_m(I)| \le    B_0 \tilde L_0 L_1^2L_2\,  C_1^{s} C_2^{m-2}
\Gamma(\rho s + (\mu-1)(m-2)  )  \, . 
\]
Now \eqref{eq:homological2} yields
\[
\begin{array}{lcrr}
\displaystyle
P_s(g_{m}(\cdot,I)) \le 
 \frac{1}{\kappa}P_{s+\tau}(B_{m}(\cdot,I))\\[0.5cm]
\displaystyle \le \, B_0e^{\tau\rho}\frac{1}{\kappa}\tilde L_0L_1^{\tau +2}L_2 C_1^{s}C_2^{m-2} 
\Gamma(\rho s + \rho\tau +(\mu-1)(m-2)  ) \, 
\end{array}
\]
for any $I\in\D^n$.  
Set $A_0:=\max\{B_0e^{\tau\rho}, 2C_0(\delta,\mu)e^{\rho}e^{\mu-1}\} $ and fix
  $$C_2\ge \frac{1}{\kappa}A_0  \tilde L_0L_1^{\tau+2}L_2.$$ 
  Since $\kappa\le 1$ and $\tilde L_0\ge 1$ the  inequality in \eqref{eq:constants1} holds as well. As $\mu= \rho(\tau+1) +1$,  we obtain 
\[
P_s(g_{m}(\cdot,I)) 
\le C_1^{s}C_2^{m-1} 
\Gamma(\rho s + (\mu-1)(m-1) -\rho )   
\]
for any $I\in\D^n$. This completes the induction and  proves Proposition \ref{prop:main-estimates}. 
 \finishproof

\noindent
{\em Proof of Proposition \ref{prop:gevrey}}. 
By the Cauchy formula and Lemma \ref{lemma:gevrey-estimates} we get for any $\alpha\in\N^n$ with $|\alpha|\ge 2$ the estimate
\[
\begin{array}{lcrr}
\displaystyle \left|\partial_\theta^\beta g_{|\alpha|,\alpha}(\theta)\right| 
=\frac{1}{\alpha !} \left|\partial_I^\alpha \partial_\theta^\beta g_{|\alpha|}(\theta,0)\right| \\[0.3cm]
\displaystyle \le \sup_{I\in \D^n} \left|\partial_\theta^\beta g_{|\alpha|}(\theta,I)\right| 
\le \sup_{I\in \D^n}P_{|\beta|}(g_{|\alpha|}(\cdot,I)).
\end{array}
\]
Now Proposition \ref{prop:main-estimates} and 
(\ref{eq:BetaFunction1}) imply 
\[
\begin{array}{lcrr}
\displaystyle \sup_{\theta\in \T^n} \left|\partial_\theta^\beta g_{|\alpha|,\alpha}(\theta)\right| \le  C_1^{|\beta|}C_2^{|\alpha|-1} 
\Gamma(\rho |\beta| + (\mu-1) (|\alpha|-1) -\rho)\\[0.3cm]
\displaystyle \le  (cC_2)^{-1}(cC_1)^{|\beta|}(cC_2)^{|\alpha|} 
\Gamma(\rho |\beta|+1)\Gamma((\mu -1)|\alpha| +1)\, ,
\end{array}
\]
for any $\alpha,\beta\in \N^n$, where $c=c(\rho,\mu)\ge 1$.  Using the Borel extension theorem in Gevrey classes  (see \cite[Theorem 3.7]{P3} for a more general version) we find a ${\mathcal G}^{\rho,\mu}$-smooth function $g$ such that $g(\theta,I)\sim \sum_{m=2}^\infty \sum_{|\alpha|=m} g_{m,\alpha} (\theta)I^\alpha$, i.e. the Taylor seres of $g$ is given by (\ref{eq:g}). Moreover, we have 
\begin{equation}
\label{eq:g-estimates}
\displaystyle \sup_{(\theta,I)\in \A'} |\partial_\theta^\beta \partial_I^\alpha g(\theta, I)| 
\displaystyle \le  \frac{C_0}{C_2} C_1^{|\beta|}C_2^{|\alpha|} 
\Gamma(\rho |\beta|+1)\Gamma(\mu |\alpha| +1)\, ,
\end{equation}
for any $\alpha,\beta\in \N^n$, where $\A'=\T^n\times D'$, $D'$ is a neighborhood of $0$ in $\R^n$,  $C_0=C_0(\rho,\tau,n)\ge 1$, and the constants $C_1\ge 1$ and $C_2\ge 1$ are equivalent to $L_1$ and $\frac{1}{\kappa} \tilde L_0L_1^{\tau+2}L_2$ respectively. 
In particular, 
$g$ belongs to ${\mathcal G}_{C_1, C_2}^{\rho,\mu}(\A')$ and 
$\|g\|_{C_1, C_2}\le C_0/C_2$.  
In the same way we find $H^0\in {\mathcal G}_{C_2}^{\mu}(D')$ such that 
$H^0(I)\sim \sum R_m (I)$.  
Then, using the  composition of Gevrey functions \cite[Proposition A.4]{P3}, we  
show that the function $ H'$ defined by $H'(\theta,I):= H(\theta, I + \partial g/\partial \theta (\theta, I))$ belongs to ${\mathcal G}_{ C_1,  C_2}^{\rho,\mu}(\A')$, where the Gevrey constants $C_1$ and $ C_2$ are equivalent to $L_1$ and to $\frac{1}{\kappa} \tilde L_0L_1^{\tau+2}L_2$ respectively. Recall that $\tilde L_0$ is equivalent to $L_0L_1^{n+2}$, hence, $C_1$ and $C_2$ satisfy \eqref{eq:gevrey-constants1}. 
This completes the proof of  Proposition  \ref{prop:gevrey}. \finishproof

\noindent
{\em Proof of Theorem  \ref{theo:gevrey}}. 
We are going to solve the equation 
\begin{equation}
\label{eq:implicit}
\varphi = \theta + \nabla_I g(\theta, I)\, ,\ (\theta,I)\in\A'\, ,
\end{equation}
with respect to $\theta \in \T^n$, by means of the implicit function theorem in anisotropic Gevrey classes \cite[Proposition A.2]{P3}. By (\ref{eq:g-estimates}) we have 
\begin{equation}
\label{eq:g-estimates1}
\displaystyle \sup_{(\theta,I)\in \T^n\times D'} \|\partial_\theta^\beta \partial_I^\alpha\nabla_I g(\theta, I)\| 
\displaystyle \le   C_0C_1^{|\beta|} C_2^{|\alpha|} 
\Gamma(\rho |\beta|+1)\Gamma(\mu |\alpha| +1)\, 
\end{equation}
where $C_0$ is equivalent to $1$, $C_1$  is equivalent to $L_1$, and $C_2$  is equivalent to $\frac{1}{\kappa}  L_0L_1^{\tau+n+4}L_2$. 
Set 
$$\epsilon:= (2C_0C_1)^{-1}<1 \quad \mbox{and}\quad \tilde C_2:= a(\rho,\mu,n)C_2/\epsilon.$$
 Then choosing  $D'$ small enough and $a(\rho,\mu,n)\gg 1$ we obtain by (\ref{eq:g-estimates1})
\[
\displaystyle \sup_{(\theta,I)\in \T^n\times D'} \|\partial_\theta^\beta \partial_I^\alpha\nabla_I g(\theta, I)\| 
\displaystyle \le  \epsilon  C_0C_1^{|\beta|} \tilde C_2^{|\alpha|} 
\Gamma(\rho |\beta|+1)\Gamma(\mu |\alpha| +1)\, 
\]
for any $\alpha,\beta\in \N^n$. 
Now, $\epsilon C_0C_1< 1/2$ and \cite[Proposition A.2]{P3} implies that there is
$ \theta\in {\mathcal G}_{C_1,  C_2}^{\rho,\mu}(\A', \T^n)$ which solves \eqref{eq:implicit}, where 
\[
C_1= c(\rho,\tau,n) L_1\ \mbox{and}\  C_2 =c(\rho,\tau,n)\frac{1}{\kappa} L_1^{\tau +n+5} L_2.
\]
Next using the  theorem of composition of Gevrey functions \cite[Proposition A.4]{P3}, we prove that 
$R^0(\varphi,I):= R(\theta(\varphi,I),I)$ belongs to the class $ {\mathcal G}_{C_1,  C_2}^{\rho,\mu}(\A', \A)$, where $C_1$ and $C_2$ are given by Remark \ref{rm:gevrey-constants}. 
By the same argument,   the canonical transformation $\chi$ generated by $g$ belongs to the class 
$ {\mathcal G}_{C_1,  C_2}^{\rho,\mu}(\A', \A)$. This  completes the proof of Theorem \ref{theo:gevrey} and of Remark \ref{rm:gevrey-constants}.

\section{Complements on the Gamma function}\label{sec:Gamma}
Here we collect 
certain estimates  of the Euler Gamma function
\begin{equation}
\Gamma (x)\ =\ \int_{0}^{\infty}\, e^{-t}\, t^{x-1}\ dt \ , x>0, 
\label{eq:GammaFunction}
\end{equation}
that have been used above. 
Recall that $\Gamma(x+1)=x\Gamma(x)$ which implies $\Gamma(m+1) = m!$ for any $m\in\N$. Moreover,  $\Gamma(t)$ is convex in the interval $(0,+\infty)$, it has a  minimum at some point $t_0\approx 1,46$ and  $\Gamma (t_0)\approx 0,89$. In particular, $\Gamma(t)$ is strictly decreasing in $(0,t_0]$ and strictly increasing in $[t_0,+\infty)$. 
We have the following relation (see \cite{kn:BE}, \cite{kn:O}) 
\begin{equation}
\label{eq:gamma-beta}
\Gamma(x)\Gamma(y)\ =\ \Gamma(x+y)B(x,y)\, ,\ x,y > 0\, ,
\end{equation}
where $B(x,y)$ is the Beta function which is  defined for $x>0$ and $y>0$ by the following integral representation 
\begin{equation}
B(x,y)\ =\ \int^1_0\, (1-t)^{x-1}\, t^{y-1}\, dt\, .
\label{eq:BetaFunction}
\end{equation}
Obviously, $B(x,y)$ is symmetric, i.e. $B(x,y)=B(y,x)$, and it is  decreasing with respect to both variables $x$ and $y$. Using the integral representation in (\ref{eq:BetaFunction}) and the Cauchy inequality  we get
for any positive numbers $a,b,c,d$ the following inequality
\begin{equation}
\label{eq:cauchy}
B(a  + b, c + d  ) 
\le  B(2a,2c  )^{1/2} B(2b,2d )^{1/2}. 
\end{equation}
Denote by  $[x]$  the entire part of $x\in \R$.
Since $B(x,y)$ is a decreasing function with respect to both variables $x,y > 0$, we have 
\begin{equation}
B(x,y) \ge B([x]+1,[y]+1)  =\frac{[x]![y]!}{([x]+[y]+1)!} \ge 4^{-x-y}.
\label{eq:BetaFunction1}
\end{equation}
For any  $x,y\ge 0$ we get in the same way 
\[
\begin{array}{lcrr}
 \displaystyle B(x+1,y+1) \le B([x]+1,[y]+1)  
= \frac{1}{[x]+[y]+1} \begin{pmatrix}[x]+[y]\cr [y]\end{pmatrix}^{-1} \\[0.5cm]
  \displaystyle <  \frac{3}{x+y+1 }\begin{pmatrix}[x]+[y]\cr [y]\end{pmatrix}^{-1}.
 \end{array}
\]
More generally, we have the following 
\begin{lemma}\label{lemma:gamma}
For any  $\nu\ge 1$ and $\delta>0$  there is a constant $C'(\nu,\delta)\ge 1$ such that
for any   $x,y\ge 0$  the following inequality holds
$$
 \displaystyle \begin{pmatrix}[x]+[y]\cr [x]\end{pmatrix}^\nu B(\nu x+\delta, \nu y+\delta)\ \leq\
\frac{C'(\nu,\delta)}{\big(\min(x+1,y+1)\big)^{(\nu+1)/2}}\,  .
$$
\end{lemma}
\noindent {\em Proof.} Since $B(x,y)$ is a decreasing function with respect to both variables $x>0$ and $y>0$ we can suppose that $\delta \le 1$.  Fix $0<\epsilon\le 1$. 
By Stirling's formula and the continuity of the Gamma function in the interval $[\epsilon,+\infty)$, there is 
$L=L(\epsilon)>1$  such that for any $t\geq \epsilon$ we
have
$$
L^{-1}\ \leq\ \Gamma(t)(2\pi)^{-1/2}t^{\frac{1}{2}-t}e^t\ \leq\ L\, .
$$
For any $\nu\ge 1$ and $t\ge \epsilon$
this implies the two-sided inequality
\begin{equation}
\label{eq:two-sided-inequality}
\displaystyle
L^{-\nu-1}\left(\frac{t}{2\pi}\right)^{(\nu-1)/2} \nu^{\nu t-\frac{1}{2}}\le \frac{\Gamma(\nu t)}{\Gamma(t)^\nu}\ \leq\
L^{\nu+1} \left(\frac{t}{2\pi}\right)^{(\nu-1)/2} \nu^{\nu t-\frac{1}{2}}\, .
\end{equation}
 Set $\epsilon := \delta/\nu\in (0,1]$. 
Substituting  $t=x+\epsilon$, $t=y+\epsilon$ and $t=x+y+2\epsilon$ in \eqref{eq:two-sided-inequality}, where $x\ge 0$ and $y\ge 0$,  we obtain
\[
\begin{array}{lcrr}
\displaystyle 
B(\nu x +\delta ,\nu y+\delta) = \frac{\Gamma(\nu(x+\epsilon))\Gamma(\nu(y+\epsilon))}{\Gamma(\nu(x+y+2\epsilon))}\\[0.5cm]
\displaystyle  \leq\ C\, 
\left(\frac{(x+\epsilon)(y+\epsilon)}{x+y+2\epsilon}\right)^{\frac{\nu-1}{2}}\, B(x+\epsilon,y+\epsilon)^\nu\\[0.5cm]
\displaystyle = C\, (x+y+1+2\epsilon)^{(\nu-1)/2}B(x+1+\epsilon,y+1+\epsilon)^{(\nu-1)/2}\\[0.5cm]
\displaystyle 
\times\,  B(x+\epsilon,y+\epsilon)^{(\nu+1)/2},
 \end{array}
\]
where $C= L^{3\nu+3}(2\pi)^{(1-\nu)/2}\nu^{-1/2}$. 
 For $y\ge x \ge 0$
we have 
\[
\begin{array}{lcrr}
\displaystyle  B(x+\epsilon,y+\epsilon) \\[0.5cm]
\displaystyle 
\le  \frac{(x+y+2\epsilon)(x+y+1+2\epsilon)}{(x+\epsilon)(y+\epsilon)}B(x+1+\epsilon,y+1+\epsilon)\\[0.5cm]
\displaystyle
\le \frac{2(x+y+1+2\epsilon)}{x+\epsilon}B(x+1+\epsilon,y+1+\epsilon)\, ,
 \end{array}
\]
which implies
\[
\begin{array}{lcrr}
\displaystyle 
B(\nu x +\delta ,\nu y+\delta) 
 \le C\,  2^{(\nu+1)/2}\, \frac{(x+y+1+2\epsilon)^{\nu}}{(x+\epsilon)^{(\nu+1)/2}}\\[0.5cm]
\displaystyle 
\times\, 
 B(x+1+\epsilon,y+1+\epsilon)^{\nu}.
 \end{array}
\]
Since $B(x,y)$ is a decreasing function with respect to both variables $x$ and $y$, we obtain 
\[
\begin{array}{lcrr}
B(x+1+\epsilon,y+1+\epsilon) \le B([x]+1,[y]+1)  \\[0.5cm]
 \displaystyle= \frac{1}{[x]+[y]+1} \begin{pmatrix}[x]+[y]\cr [y]\end{pmatrix}^{-1} 
 \le  \frac{5}{x+y+1+ 2\epsilon }\begin{pmatrix}[x]+[y]\cr [y]\end{pmatrix}^{-1}.
 \end{array}
\]
This implies 
\[
\displaystyle 
B(\nu x +\delta ,\nu y+\delta) 
 \le  \, 
\frac{C'}{(x+1)^{(\nu+1)/2}}\begin{pmatrix}[x]+[y]\cr [x]\end{pmatrix}^{-\nu}.
\]
where $C'= \left(\frac{2\nu}{\delta}\right)^{(\nu+1)/2}5^\nu C$. This 
 completes the  proof of the assertion since the inequality in Lemma \ref{lemma:gamma} is symmetric with respect to $x,y$. \finishproof
\begin{rem}\label{rem:gamma}
As in  \eqref{eq:two-sided-inequality} one  proves that for any $\rho>0$ there is a constant $C(\rho)>1$ such that
\[
C(\rho)^{-m}\Gamma(\rho m+1)\le m !\,  ^\rho  \le C(\rho)^{m}\Gamma(\rho m+1)
\]
 for any $m\in\N$. 
\end{rem}

We have also
\begin{lemma}\label{lemma:binom} 
For any $m\in \N$ and $\alpha\in \N^n$ such that $2\le |\alpha|\le  m$  we have 
\[
\displaystyle
 \sum_{(\alpha^1,\ldots ,\alpha^{m})\in 
 {\widetilde\N}(\alpha,m )} \frac{\alpha!}{\alpha^1 ! \ldots \alpha^{m}! } 
= \frac{(m-1)!}{(m-|\alpha |)!(|\alpha|-1)!}\,  ,
\]
where ${\N}(\alpha,m )$ is defined by \eqref{eq:index-set}.
\end{lemma}

\noindent
{\em Proof}. Let $f$ be 
an analytic function  in ${\R}^{n}$.
Then the left hand-side of the inequality above coincides with
the coefficient $a_{\alpha,m} $ in the identity
\[
\begin{array}{lcrr}
\displaystyle f\left( \frac{X}{1-X}, \ldots ,\frac{X}{1-X}\right) 
 = \sum_{\gamma\in \N^{n} }\left( \sum_{k=1}^\infty X^k,\ldots, \sum_{k=1}^\infty X^k\right)^\gamma
\frac{f^{(\gamma)}(0)}{\gamma!}\\[0.5cm]
\displaystyle =
\sum_{\gamma\in
\N^{n} }\sum_{\nu\ge |\gamma|}^{\infty} a_{\gamma,\nu}\, X^\nu\, 
\frac{f^{(\gamma)}(0)}{\gamma!} \, .
 \end{array}
\]
Now taking $f(Y)=Y^\alpha$ we obtain 
$$
a_{\alpha,m}  = \frac{1}{m !}\left( \frac{d}{dX}
\right)^{m} \left( X^{|\alpha|}(1-X)^{-|\alpha| }\right)
\mid_{X=0} = \frac{(m -1)!}{(m -|\alpha|)! (|\alpha|-1)!}\,  .
$$
\finishproof

\noindent
\textbf{Acknowledgments.}
{\em Part of this work has been done at the Institute of Mathematics, Bulgarian Academy of Sciences, and I would like to thank the colleagues there for the stimulating discussions. }

\vspace{0.5cm} 
\noindent 
T.M.:
University of Rousse, \\
Department of Algebra and Geometry,\\
7012 Rousse, Bulgaria

\vspace{0.5cm} 
\noindent 
G. P.: 
Universit\'e de Nantes,  \\
Laboratoire de mathématiques Jean Leray,\\
2, rue de la Houssini\`ere,\\
 BP 92208,  44072 Nantes 
Cedex 03, France \\
e-mail: georgi.popov@univ-nantes.fr\\
and\\
Institute of Mathematics and Informatics\\
Bulgarian Academy of Sciences\\
"Acad. G. Bonchev" Str. 8\\
1113 Sofia, Bulgaria
 
\vspace{0.5cm} 
\noindent 
\end{document}